\documentclass[nameyear,preprint,review,twocolumn]{autart} 

\usepackage{bbm,amsmath,mathtools,amssymb,graphicx}          
\usepackage{natbib}

\begin{document}

\begin{frontmatter}

\title{Output Consensus of Networked Hammerstein and Wiener Systems\thanksref{footnoteinfo}} 

\thanks[footnoteinfo]{The work was supported by the National Key Basic Research Program of China (973 program, 2014CB845301), and the National Center for Mathematics and Interdisciplinary Sciences, Chinese Academy of Sciences. Corresponding author Wenhui Feng Tel. +86-15210965635}

\author[feng]{Wenhui Feng}\ead{feng\_wh@amss.ac.cn},    
\author[chen]{Han-Fu Chen}\ead{hfchen@iss.ac.cn},               

\address[feng]{Key Laboratory of Systems and Control, Institute of System Science, Academy of Mathematics and Systems Science, Chinese Academy of Sciences, University of Chinese Academy of Sciences, 55 Zhongguancundonglu, Beijing 100190, PR China}  
\address[chen]{Key Laboratory of Systems and Control, Institute of System Science, Academy of Mathematics and Systems Science, Chinese Academy of Sciences, 55 Zhongguancundonglu, Beijing 100190, PR China}             


\begin{abstract}                          
In this paper we consider the output consensus problem of networked Hammerstein and Wiener systems in a noisy environment. The Hammerstein or Wiener system is assumed to be open-loop stable, and its static nonlinearity is allowed to grow up but not faster than a polynomial. A control algorithm based on the distributed stochastic approximation algorithm with expanding truncations is designed and it is shown that under the designed control the output consensus is achieved. The numerical simulation given in the paper justifies the theoretical assertions.
\end{abstract}

\end{frontmatter}

\section{Introduction}
In past decades the consensus problem for multi-agent systems has drawn much attention from researchers for its close connection with problems arising from biological science, physical science, computer science, and other areas. Among the early theoretical studies the work \cite{Jadbabaie2003} gives the mathematical explanation for the physical
phenomenon discovered by \cite{Vicsek1995}, pointing out that the problem can be reduced to stability analysis for a class of first order integrator systems. This problem is then considered in \cite{Saber2004} for cases including the fixed directed graph, the switched directed graphs, and the undirected graph with time delay. Besides, the concept of algebraic connectivity is expanded in \cite{Saber2004} from the undirected graph to the balanced directed graph, which plays an important role in achieving averaged consensus. More general adjacency matrices and Laplace matrices in connection with the consensus of multi-agents are considered in \cite{Ren2005}, where it is pointed out that for time varying graphs the asymptotic consensus can be guaranteed if the union graph contains a spanning tree. All works mentioned above are mainly with the first-order linear systems in a noise-free communication environment.

Since the noise is unavoidable in practice, the noise environment has naturally been taken into account in later research. The fixed directed graph with white observation noise is dealt with in \cite{Li2009}, where the sufficient and necessary conditions are obtained to guarantee the asymptotic unbiased mean square consensus. Further, for the time varying graph the sufficient conditions are given in \cite{Li2010} to ensure the mean square and almost sure consensus. The almost sure consensus for the first-order discrete-time systems on fixed graph with observation noise is transformed to the convergence analysis of some stochastic approximation algorithm in \cite{Huang2009}. The noise conditions used in \cite{Huang2009} and \cite{Li2010} have been weakened in \cite{Fang2012}. Since the second-order systems are of explicit physical meaning, consensus of the second-order multi-agent systems has also attracted attention from many researchers \cite{ChenY2013}, \cite{Ren2008}, and \cite{Yu2017}.

The works concerned above are with linear systems, but in practice nonlinear systems are ubiquitous. The consensus problem for some nonlinear systems can be reduced to the consensus of linear multi-agent systems by local linearization of the nonlinearity if it satisfies the Lipschitz condition (\cite{Saber2007}). The linearization approach is hard to work when the Lipschitz condition does not take place. The consensus problem of multi-agent nonlinear systems is considered in many papers e.g., \cite{Hua2016}, \cite{Li2013}, \cite{Liu2013}, \cite{Liu2017}, \cite{Munz2011}, and \cite{Wang2017}, where various types of consensus including the leader-following consensus, finite-time consensus, adaptive consensus etc. are discussed under different settings. For example, the tracking consensus for a class of high-order nonlinear systems with unknown parameters and external disturbances is considered in \cite{Wang2017}, where the distributed adaptive control based on back-stepping method is given so that the boundedness of the closed-loop system and the output tracking consensus are achieved.

In general, for the consensus problem of nonlinear systems to guarantee stability of the closed-loop systems is of primary importance. For this the growth rate of system nonlinearity, the noise, and the step-size used in the algorithm play an important role. When the classical Lyapunov method is used in stability analysis, rather strong conditions are usually required. For example, in the leader-following case a stable leader is needed to serve as a reference signal in order to avoid divergence of agents, and thus the problem turns to tracking consensus for multi-agents (\cite{Hua2016}, \cite{Liu2017}, and \cite{Wang2017}), while in other cases the stability of closed-loop systems is normally guaranteed by imposing the Lipschitz condition on nonlinearity (\cite{Li2013} and \cite{Liu2013}).

In this paper the nonlinear dynamical systems of agents are either the Hammerstein system or the Wiener system, which have the wide background in practice. In contrast to the conditions discussed above, here neither a stable leader nor the Lipschitz condition are needed. As mentioned in \cite{Fang2012}, \cite{Huang2009}, and \cite{Li2010}, there is a close relationship between the consensus problem of multi-agents and the root-seeking problem for an unknown function. We apply the distributed stochastic approximation algorithm with expanding truncations (DSAAWET) to solve the output consensus problem for networked Hammerstein and Wiener systems. DSAAWET used in the paper is not exactly the same as that proposed in \cite{Lei2016} and \cite{Lei2015}, but the basic idea remains the same.

The rest of the paper is organized as follows. The problem for output consensus of networked Hammerstein and Wiener systems is described in Section 2. The control is defined by a DSAAWET in Section 3. The properties of noises appearing in DSAAWET are analyzed in Section 4. In Section 5, the auxiliary sequences are introduced. Convergence of the algorithm is proved in Section 6. Simulation and conclusions are given in Sections 7 and 8, respectively.

\subsection*{Notations}
Let $\mathbb{R}$ denote the real line, and $\mathbb{R}^{n}$ the linear space of $n$-dimensional vectors. $\mathbbm{1}_{n}$ denotes the $n$-dimensional vector with all entries equal to one. For $\mathbbm{x}=[x_{1},\cdots,x_{N}]^T\in \mathbb{R}^N$, define the norm $\|\mathbbm{x}\|_{\infty}=\max_{1\le i\le N}\{|x_{i}|\}$. By $\mathbb{R}^{m\times n}$ we denote the linear space of $m\times n$ matrices, and by $\text{Null}(A)$ the null space of matrix $A$. Let $[s]$ denote the integer part of a nonnegative real number $s$, and $\text{Span}\{\mathbbm{e}_{1},\cdots,\mathbbm{e}_{n}\}$ be the linear subspace spanned by vectors $\mathbbm{e}_{1},\cdots,\mathbbm{e}_{n}$. Set $p\wedge q\triangleq\min\{p,q\}$, and $p\vee q\triangleq\max\{p,q\}$.

\section{Problem Description}
Consider a network of $N$ agents and the corresponding topology being an undirected graph $\mathcal{G}=(\mathcal{N},\mathcal{E})$, where $\mathcal{N}=\{ 1,2,\cdots,N \}$ is the node (agent) set and $\mathcal{E} \subset \mathcal{N} \times \mathcal{N}$ is the edge set. The neighbor set of agent $i$ is denoted by $\mathcal{N}_{i}=\{ j|(j,i) \in \mathcal{E} \}$ and it is assumed that $i \notin \mathcal{N}_{i}$, i.e., the graph $\mathcal{G}$ does not contain any self-loops. A path from agent $j$ to agent $i$ is denoted by $j=i_{0},i_{1},\cdots,i_{d_{ij}-1},i_{d_{ij}}=i$, where $(i_{s},i_{s+1}) \in \mathcal{E}, s=0,1,\cdots,d_{ij}-1$ and $d_{ij}$ is called the length of this path. The shortest length of these paths is called the distance from $j$ to $i$, still denoted as $d_{ij}$. $d(\mathcal{G})\triangleq\max\{d_{i,j}, i,j \in \mathcal{N}\}$ is the diameter of $\mathcal{G}$. $P=[p_{ij}]_{N \times N}$ is called the adjacency matrix, if $p_{ij}>0, j \in \mathcal{N}_{i}; p_{ij}=0, j \notin \mathcal{N}_{i}$. Set $p_{i}\triangleq\sum_{j \in \mathcal{N}_{i}}p_{i,j}$ and $D=\text{diag}\{p_{1}, \cdots ,p_{N}\}$ which is called the degree matrix. $L\triangleq D-P$ is the Laplace matrix of $\mathcal{G}$.

The dynamics of agent $i \in \mathcal{N}$ is the following SISO discrete-time Hammerstein or Wiener system.

Hammerstein system:
\begin{equation}\label{HSystem}
v_{i,k}=f_{i}(u_{i,k}),\text{ }C_{i}(z)y_{i,k+1}=D_{i}(z)v_{i,k}
\end{equation}

Wiener system:
\begin{align}\label{WSystem}
C_{i}(z)v_{i,k+1}=D_{i}(z)u_{i,k},\text{ }y_{i,k+1}=f_{i}(v_{i,k+1}),
\end{align}
where $u_{i,k}, v_{i,k}, y_{i,k} \in \mathbb{R}$ are the control input, internal variable, and output, respectively. The first subscript $i$ represents agent and the second subscript $k$ represents the discrete-time. $f_{i}(\cdot):\mathbb{R} \to \mathbb{R}$ is an unknown function. By $z$ we denote the backward shift operator: $zy_{i,k+1}=y_{i,k}$. In \eqref{HSystem} and \eqref{WSystem}
\[
C_{i}(z)=1+c_{i,1}z+ \cdots +c_{i,p}z^p
\]
and
\[
D_{i}(z)=1+d_{i,1}z+ \cdots +d_{i,q}z^{q}
\]
are polynomials of $z$, where $c_{i,s}, d_{i,r}\in \mathbb{R},s=1, 2, \cdots ,p, r=1, 2, \cdots ,q$ are unknown parameters, and $p$ and $q$ are also unknown and may depend on $i$ but we still denote them as $p$ and $q$ for simplicity of writing.

The observation of neighbor $j \in \mathcal{N}_{i}$ at agent $i \in \mathcal{N}$ is
\begin{equation}\label{Observations}
z_{ij,k+1}=y_{j,k+1}+\epsilon_{ij,k+1},
\end{equation}
where $\epsilon_{ij,k+1}$ is the observation noise.

Define the function
\begin{align*}
h_{i}(u)=\begin{cases}\frac{d_{i}}{c_{i}}f_{i}(u),&\text{if agent }i\text{ is the Hammerstein system,}\\
f_{i}(\frac{d_{i}}{c_{i}}u),&\text{if agent }i\text{ is the Wiener system,}\end{cases}
\end{align*}
where $c_{i} \triangleq 1+\sum_{s=1}^p c_{i,s},\text{ } d_{i} \triangleq 1+\sum_{s=1}^{q} d_{i,s}$. By the following condition A1, $c_{i}=C_{i}(1) \ne 0$.

We list the conditions to be used.
\begin{itemize}
\item[A1:] $C_{i}(z), \text{ }i \in \mathcal{N}$ are stable, i.e., the roots of $C_{i}(z)$ are outside the unit disk.
\item[A2:] \begin{itemize}
            \item[i)] $f_{i}(\cdot),\text{ }i \in \mathcal{N}$ are continuous.
            \item[ii)] There exists an unknown constant $\mu > 0$ such that
                \[
                | f_{i}(u) |=O(| u |^\mu) \text{ as } |u| \to \infty \text{ }\forall i \in \mathcal{N}.
                \]
            \item[iii)] $h_{i}(u):\mathbb{R} \to \mathbb{R}, i \in \mathcal{N}$ are strictly monotonically increasing and have range $(-\infty,+\infty)$.
          \end{itemize}
\item[A3:] $\{ \epsilon_{ij,k} \}_{k\ge1}, i\in \mathcal{N}, j\in\mathcal{N}_{i}$ are mutually independent sequences with
    \begin{align}\label{ConditionForNoise}
    \mathbb{E}\epsilon_{ij,k}=0,\text{ }\sup_{k}\mathbb{E}|\epsilon_{ij,k}|^2 < \infty.
    \end{align}
\item[A4:] $\mathcal{G}$ is a connected and undirected graph.
\end{itemize}

The output consensus problem is stated as follows: at agent $i \in \mathcal{N}$ the control input $u_{i,k+1}$ should be designed on the basis of its output $y_{i,k+1}$ and the observations on neighbors $z_{ij,k+1},j \in \mathcal{N}_{i}$ so that the outputs of all agents converge to the same limit
\begin{equation}\label{ConsensusTarget}
\lim_{k \to \infty} \mathbbm{y}_{k}=y^{0} \mathbbm{1}_{N} \text{ a.s.},
\end{equation}
where $\mathbbm{y}_{k}=[y_{1,k},\cdots,y_{N,k}]^T$ and $y^{0} \in \mathbb{R}$ may depend on samples $\omega \in \Omega$.

Under A1 and A2 i) by Lemma~1 in \cite{Chen2007} it is known that for both \eqref{HSystem} and \eqref{WSystem}, $\lim_{k\to\infty}y_{i,k}=h_{i}(u_{i})$ as $\lim_{k\to\infty}u_{i,k}=u_{i}$.

Set
\[
\mathbbm{u}_{k}=[u_{1,k},\cdots,u_{N,k}]^{T}
\]
and
\[
\mathbbm{h}(\mathbbm{u}_{k})=[h_{1}(u_{1,k}),\cdots,h_{N}(u_{N,k})]^{T}.
\]
Then, the problem is restated as follows: at agent $i \in \mathcal{N}$ the control input $u_{i,k+1}$ should be designed on the basis of its output $y_{i,k+1}$ and the observations $z_{ij,k+1},j \in \mathcal{N}_{i}$ on neighbors so that
\begin{align}\label{LimUk}
\lim_{k\to\infty}\mathbbm{u}_{k}=\mathbbm{u}^{0}\in J\triangleq\{\mathbbm{u}\in\mathbb{R}^{N}: \mathbbm{h}(\mathbbm{u})\in\text{Span}\{\mathbbm{1}_{N}\}\},
\end{align}
where $\mathbbm{u}^{0}\triangleq[u_{1}^{0}, \cdots, u_{N}^{0}]^{T}$ may depend on samples.

For solvability of the consensus problem, $J$ in \eqref{LimUk} must be nonempty, i.e., the intersection of ranges of $h_{i}(\cdot), i=1,\cdots,N$ should be nonempty. Besides, the steady output may depend on samples since communication between agents is corrupted by noises. In order to reach consensus for almost all samples, the strong condition A2 iii) is imposed on $h_{i}(\cdot)$.

Set $c_{i,s}=d_{i,r}=0$ for $s>p$ and $r>q$ and define matrices
\[
C_{i} \triangleq  \begin{bmatrix}
                       -c_{i,1} & 1     & 0      & \cdots & 0          \\
                       \vdots   &       & \ddots & \ddots & \vdots   \\
                       \vdots   &       &        & \ddots & 0        \\
                       \vdots   &       &        &        & 1        \\
                       -c_{i,(p \vee q)+1} & 0     & \cdots & \cdots & 0
               \end{bmatrix},
\]
\[
D_{i} \triangleq \begin{bmatrix}
                     1 & d_{i,1} & \cdots & d_{i,p \vee q}
              \end{bmatrix}^T,
\]
\[
G_{1} \triangleq [1  \hspace{2.5mm} 0 \cdots 0]_{1\times ((p \vee q) +1)}.
\]

The dynamical systems of \eqref{HSystem} and \eqref{WSystem} can be written in the state space form
\begin{align*}
Y_{i,k+1}&=C_{i}Y_{i,k}+D_{i}v_{i,k},\text{ }y_{i,k}=G_{1}Y_{i,k},\\
V_{i,k+1}&=C_{i}V_{i,k}+D_{i}u_{i,k}, \text{ }v_{i,k}=G_{1}V_{i,k}.
\end{align*}
Then it follows that
\begin{align}
Y_{i,s+l+1}&=C_{i}^{l+1}Y_{i,s}+\sum_{k=s}^{s+l}C_{i}^{s+l-k}D_{i}v_{i,k},\text{ }\forall s,\text{ }l \ge 0, \label{ExpansionForLSSY}\\
V_{i,s+l+1}&=C_{i}^{l+1}V_{i,s}+\sum_{k=s}^{s+l}C_{i}^{s+l-k}D_{i}u_{i,k},\text{ }\forall s,\text{ }l \ge 0. \label{ExpansionForLSSV}
\end{align}
It is well known that under A1, $C_{i}$ are stable and there are $r>0$ and $\delta>0$ such that
\begin{align}\label{EstimateforC}
\| C_{i}^k \| \le re^{- \delta k}\text{ }\forall k \ge 0\text{ }\forall i \in \mathcal{N}.
\end{align}
For Hammerstein systems \eqref{HSystem} from \eqref{ExpansionForLSSY} and \eqref{EstimateforC} it follows that
\begin{align}
\|Y_{i,k}\|&\le \|C_{i}^{k}\|\|Y_{i,0}\|+\sum_{s=0}^{k-1}\|C_{i}^{k-s-1}\|\|D_{i}\||v_{i,s}|\notag\\
&\le re^{-\delta k}\|Y_{i,0}\|+r\sup_{0\le s\le k-1}|v_{i,s}|\|D_{i}\|. \label{NorOfYik}
\end{align}
Similarly, for Wiener systems \eqref{WSystem}, from \eqref{ExpansionForLSSV} and \eqref{EstimateforC} we have
\begin{align}
\|V_{i,k}\|&\le \|C_{i}^{k}\|\|V_{i,0}\|+\sum_{s=0}^{k-1}\|C_{i}^{k-s-1}\|\|D_{i}\||u_{i,s}|\notag\\
&\le re^{-\delta k}\|V_{i,0}\|+r\sup_{0\le s\le k-1}|u_{i,s}|\|D_{i}\|. \label{NorOfVik}
\end{align}

\section{Control algorithm}
Define $\mathbbm{g}(\mathbbm{u})\triangleq[g_{1}(\mathbbm{u}),\cdots,g_{N}(\mathbbm{u})]^{T}:\mathbb{R}^{N} \to\mathbb{R}^{N}$, where
\begin{align*}
g_{i}(\mathbbm{u})&\triangleq\sum_{j\in\mathcal{N}_{i}}p_{ij}(h_{j}(u_{j})-h_{i}(u_{i})) \\
& =\sum_{j\in\mathcal{N}_{i}}p_{ij}h_{j}(u_{j})-p_{i}h_{i}(u_{i}).
\end{align*}

By A4, it is noted in \cite{Saber2004} and \cite{Ren2005} that $L=D-P\in\mathbb{R}^{N\times N}$ is nonnegative definite and
\begin{align}\label{NullLEqSpan}
\text{Null}\{L\}=\text{Span}\{\mathbbm{1}_{N}\}.
\end{align}
Since $\mathbbm{g}(\mathbbm{u})=-L\mathbbm{h}(\mathbbm{u})$, by \eqref{NullLEqSpan} it follows that
\begin{align}\label{RootsOfG}
\mathbbm{g}(\mathbbm{u})=0 \Leftrightarrow \mathbbm{h}(\mathbbm{u})\in\text{span}\{\mathbbm{1}_{N}\} \Leftrightarrow \mathbbm{u}\in J,
\end{align}
i.e., the root set of $\mathbbm{g}(\mathbbm{u})$ coincides with $J$ defined in \eqref{LimUk}. Therefore, the output consensus problem is transformed to root-seeking for some function. The control algorithm will be constructed on the basis of DSAAWET. The regression function of agent $i$ is $g_{i}(\mathbbm{u})$.

Assume the estimate for the roots of $\mathbbm{g}(\mathbbm{u})$ at time $k$ is $\mathbbm{u}_{k}=[u_{1,k},\cdots,u_{N,k}]^{T}$. Then, as observation of $g_{i}(\mathbbm{u}_{k})$ at $i$ we may take
\begin{align}\label{ObservationAgenti}
O_{i,k+1}&=\sum_{j\in\mathcal{N}_{i}}p_{ij}(z_{ij,k+1}-y_{i,k+1})\notag\\ &=\sum_{j\in\mathcal{N}_{i}}p_{ij}y_{j,k+1}-p_{i}y_{i,k+1} +\sum_{j\in\mathcal{N}_{i}}p_{ij}\epsilon_{ij,k+1}\notag \\
&=g_{i}(\mathbbm{u}_{k})+\epsilon_{i,k+1},
\end{align}
where $\epsilon_{i,k+1}=\epsilon_{i,k+1}^{(1)}+\epsilon_{i,k+1}^{(2)} +\epsilon_{i,k+1}^{(3)}$ is the observation noise, where
\begin{align}
\epsilon_{i,k+1}^{(1)}&=\sum_{j\in\mathcal{N}_{i}}p_{ij}\epsilon_{ij,k+1},\label{Noise1}\\ \epsilon_{i,k+1}^{(2)}&=p_{i}(h_{i}(u_{i,k})-y_{i,k+1}),\label{Noise2}
\end{align}
and
\begin{align}
\epsilon_{i,k+1}^{(3)}&=\sum_{j\in\mathcal{N}_{i}}p_{ij}(y_{j,k+1}-h_{j}(u_{j,k}))\label{Noise3}. \end{align}

Similar to DSAAWET proposed in \cite{Lei2016}, \cite{Lei2015}, and \cite{Fang2001}, we construct the distributed root-seeking algorithm:
\begin{align}
\sigma_{i,k}^{'}&=\max \{ \sigma_{j,k},j \in \mathcal{N}_{i},\sigma_{i,k} \}, \sigma_{i,1}=0, k\ge1, \label{SigmakHat}\\
u_{i,k}^{'}&=u_{i,k}I_{\{\sigma_{i,k}^{'}=\sigma_{i,k}\}} +u_{i}^{*}I_{\{\sigma_{i,k}^{'}>\sigma_{i,k}\}},\label{UyiPie} \\
O_{i,k+1} & =\sum_{j \in \mathcal{N}_{i}}p_{ij}(z_{ij,k+1}-y_{i,k+1}) , \label{ObservationOfi}\\
u_{i,k+1} & =\left(u_{i,k}^{'}+a_{k}O_{i,k+1}\right)I_{\{ | u_{i,k}^{'}+a_{k}O_{i,k+1} | < M_{\sigma_{i,k}^{'}} \}} \notag \\
&+ u_{i}^*I_{\{ | u_{i,k}^{'}+a_{k}O_{i,k+1} | \ge M_{\sigma_{i,k}^{'}} \}},\label{AlgorithmUk}\\
\sigma_{i,k+1} & = \sigma_{i,k}^{'} + I_{\{ | u_{i,k}^{'}+a_{k}O_{i,k+1} | \ge M_{\sigma_{i,k}^{'}} \}},\label{AlgorithmSigmak}
\end{align}
where $I_{ \{ \cdot \} }$ is the indicator function, $a_{k}=\frac1k, M_{k}=\ln (k+c_{M})$, and $c_{M}>0$ is a constant such that $|u_{i}^{*}|<M_{0},i=1,2,\cdots,N$. The new information obtained by agent $i$ at $k+1$ is $\{ z_{ij,k+1}, \sigma_{j,k}, j \in \mathcal{N}_{i} \}$ and the output $y_{i,k+1}$. The algorithm is distributed, since at each agent only the local information is used.

We explain the algorithm. The information obtained by agent $i \in \mathcal{N}$ at $k+1$ includes the observations $z_{ij,k+1}$ and truncation numbers $\sigma_{j,k}$ for neighbors in addition to itself's output and truncation number. In order to ensure control algorithms of agents update gradually and simultaneously, we need the differences of truncation numbers at agents are not too large. In the ideal case the truncations at all agents occur at the same time. However, the truncation numbers at agents cannot be always the same, since only the local information can be used and no global information is available. To make sure the differences of truncation numbers among agents are not too big we first introduce $\sigma_{i,k}^{'}=\max \{ \sigma_{j,k},j \in \mathcal{N}_{i},\sigma_{i,k} \}$ and then force the truncation number of agent $i$ at $k+1$ to catch up with its neighbors at $k$. In other words, $\sigma_{i,k+1}$ is set to equal to the largest one of $\sigma_{j,k},j\in\mathcal{N}_{i}$ when agent $i$ finds the truncation number $\sigma_{i,k}$ is smaller than $\sigma_{i,k}^{'}$, even if the estimate is still within the truncation bound. This is what \eqref{AlgorithmSigmak} means.

\section{Properties of noises}
Noticing the choice of step-size and truncation bound, similar to Lemma~2 in \cite{Chen2007} we have the following lemma.
\begin{lem}\label{BasicProperty}
Assume A1,A2 i), A2 ii), and A3. For any $l:l=0,1,2,\cdots, [\ln s]+m$ with $m$ being a given positive integer, the following limits take place:
\begin{align}
\sum_{k=s}^{s+l}a_{k} & \xrightarrow[s \to \infty]{} 0,\label{AkToLn}\\
| \sum_{k=s}^{s+l}a_{k}h_{i}(u_{i,k}) | & \xrightarrow[s \to \infty]{} 0 \text{ } \forall i\in\mathcal{N},\label{AkVkToLn}\\
| \sum_{k=s}^{s+l}a_{k}y_{i,k+1} | & \xrightarrow[s \to \infty]{} 0 \text{ } \forall i\in\mathcal{N},\label{AkYkToLn}\\
| \sum_{k=s}^{s+l}a_{k}g_{i}(\mathbbm{u}_{k}) | & \xrightarrow[s \to \infty]{} 0 \text{ }\forall i\in\mathcal{N},\label{AkGkToLn}\\
| \sum_{k=s}^{s+l}a_{k}\epsilon_{i,k+1} | & \xrightarrow[s \to \infty]{} 0 \text{ a.s.}\text{ } \forall i\in\mathcal{N},\label{AkEPkToLn}\\
| \sum_{k=s}^{s+l}a_{k}O_{i,k+1} | & \xrightarrow[s \to \infty]{} 0 \text{ a.s.} \text{ }\forall i\in\mathcal{N}.\label{AkOkToLn}
\end{align}
\end{lem}

\begin{pf}
The proof of \eqref{AkToLn} follows from the following chain of inequalities and equality:
\begin{align*}
\sum_{k=s}^{s+l}a_{k} &\le \sum_{k=s}^{s+[\ln s]+m} \frac 1k \le \sum_{k=s}^{s+[\ln s]+m} \int_{k-1}^{k} \frac {\mathrm{d}x}{x} \notag \\
&= \ln \frac {s+[\ln s]+m}{s-1} \xrightarrow[s \to \infty]{} 0.
\end{align*}

Noticing $\sigma_{i,k} \le k-1 \text{ } \forall i\in\mathcal{N}$, we have
\begin{align}\label{EstUik}
|u_{i,k}|\le M_{\sigma_{i,k}^{'}}\le M_{k-1}=\ln(k-1+c_{M}),
\end{align}
which by A2 ii) means that there exists $\alpha_{1}>0$ such that
\begin{align}\label{EstForVk}
|h_{i}(u_{i,k})| \le \alpha_{1} (\ln k)^{\mu}.
\end{align}
As a result, for $v=0,1,\cdots,q$ and $l=0,1,\cdots,[\ln k] +m$ we have
\begin{align}\label{PrForAkVk}
&|\sum_{k=s}^{s+l}a_{k}h_{i}(u_{i,k})| \le \alpha_{1} \sum_{k=s}^{s+[\ln{s}]+m} \frac {(\ln{k})^{\mu}}{k} \notag \\ &\le \alpha_{1} \sum_{k=s}^{s+[\ln{s}]+m} \int_{k-1}^{k} \frac {(\ln{x})^{\mu}}{x} \mathrm{d}x \notag \\ &= \int_{s-1}^{s+[\ln{s}]+m} (\ln{x})^{\mu} \mathrm{d}(\ln{x}) \notag \\ &= \alpha_{1} \frac{1}{1+\mu} ( (\ln{(s+[\ln{s}]+N)})^{1+\mu} -(\ln{(s-1)})^{1+\mu}) \notag \\ &= \frac{\alpha_{1}}{1+\mu} (\ln{(s-1)})^{1+\mu} ( (1 + \frac{\ln{(1+ \frac{[\ln{s}]+m+1}{s-1}})}{\ln{(s-1)}})^{1+\mu} -1) \notag \\ &= (\ln{(s-1)})^{1+\mu} O(\frac{\ln{(1+\frac{[\ln{s}]+m+1}{s-1}})}{\ln{(s-1)}}) \notag \\ &=(\ln{(s-1)})^{\mu}\cdot O(\frac{[\ln{s}]+m+1}{s-1})\notag \\ & = O(\frac{(\ln{s})^{1+\mu}+(m+1)(\ln{s})^{\mu}}{s-1})\xrightarrow[s \to \infty]{} 0,
\end{align}
and thus \eqref{AkVkToLn} is proved.

In view of A2 ii), \eqref{NorOfYik}, \eqref{NorOfVik}, and \eqref{EstUik}, in both cases \eqref{HSystem} and \eqref{WSystem} there exists $\alpha_{2}>0$ such that
\begin{align}\label{EstForYk}
|y_{i,k+1}| < \alpha_{2}(\ln k)^{\mu}
\end{align}
for sufficiently large $k$. Then similar to \eqref{PrForAkVk}, \eqref{AkYkToLn} can be proved.

By the definition of $g_{i}(\mathbbm{u}_{k})$, there exists $\alpha_{3}>0$ such that
\[
|g_{i}(\mathbbm{u}_{k})| \le \alpha_{3}(\ln k)^{\mu}
\]
for sufficiently large $k$. So, the proof for \eqref{AkGkToLn} can also be carried out in a way similar to that for \eqref{PrForAkVk}.

From A3 and $a_{k}=\frac{1}{k}$ it is clear that
\begin{align}\label{SumAkEpijk}
|\sum_{k=1}^{\infty}a_{k}\epsilon_{ij,k+1}| < \infty \text{ a.s. }\forall i,j \in \mathcal{N},
\end{align}
which implies
\begin{align*}
|\sum_{k=1}^{\infty}a_{k}\epsilon_{i,k+1}^{(1)}|&=| \sum_{k=1}^{\infty}a_{k}\sum_{j\in \mathcal{N}_{i}}p_{ij}\epsilon_{ij,k+1}|< \infty \text{ a.s.}
\end{align*}
In view of \eqref{AkVkToLn} and \eqref{AkYkToLn}, we have $|\sum_{k=s}^{s+l}a_{k}\epsilon_{i,k+1}^{(l)}| \xrightarrow[s \to \infty]{}0,l=2,3$. Then \eqref{AkEPkToLn} follows from the definition of $\epsilon_{i,k}$.

Finally, \eqref{AkOkToLn} is straightforwardly obtained from  \eqref{AkGkToLn} and \eqref{AkEPkToLn}.\qed
\end{pf}

As in \cite{Lei2016} and \cite{Lei2015}, for a positive integer $m\ge0$, set $r(i,m)=\inf \{n\ge1,\sigma_{i,n}\ge m \}$, the smallest time for agent $i$'s truncation number to reach $m$ and $r(m)=\inf \{r(i,m), i \in \mathcal{N}\}$, the smallest time for some agent's truncation number to reach $m$. So $r(i,0)=r(0)=1$. Set $\inf\emptyset=\infty$.

We now prove that
\begin{align}\label{RimRmLeDG}
0\le r(i,m)-r(m) \le d(\mathcal{G}) \text{ }\forall i \in \mathcal{N} \text{ when }r(m) < \infty,
\end{align}
where $d(\mathcal{G})$ is the diameter of $\mathcal{G}$. The first inequality of \eqref{RimRmLeDG} is seen by the definitions of $r(i,m)$ and $r(m)$. Assume $r(j,m)=r(m)$, i.e., the truncation number of agent $j$ first reaches $m$. The shortest path from $j$ to $i$ is denoted as $j=i_{0},i_{1},\cdots,i_{d_{ij}}=i$ where $(i_{s},i_{s+1})\in\mathcal{E}, s=0,1,\cdots,d_{ij}$. From
\eqref{SigmakHat}, $\sigma_{i_{s+1},r(i_{s},m)}^{'} \ge \sigma_{i_{s},r(i_{s},m)}=m$ and by
\eqref{AlgorithmSigmak}, $\sigma_{i_{s+1},r(i_{s},m)+1} \ge \sigma_{i_{s+1},r(i_{s},m)}^{'}$ which implies that $\sigma_{i_{s+1},r(i_{s},m)+1} \ge m$, i.e., the truncation number of agent $i_{s+1}$ is bigger than or equal to $m$ at $r(i_{s},m)+1$. Therefore,
\begin{align}\label{Ris1mRism}
r(i_{s+1},m) \le r(i_{s},m)+1, \text{ }s=0,1,\cdots,d_{ij}-1.
\end{align}
Noticing the path from $j$ to $i$, by \eqref{Ris1mRism} we have $r(i,m)=r(i_{d_{ij}},m) \le r(i_{d_{ij}-1},m)+1 \le \cdots \le r(i_{0},m)+d_{ij}=r(m)+d_{ij} \le r(m)+d(\mathcal{G})$, which implies the second inequality of \eqref{RimRmLeDG}.

Define the positive integer $m(k,T)\triangleq\max\{m:\sum_{s=k}^{m}a_{s} \le T\}$ for a given $T>0$. It directly follows from definition that $\sum_{s=k}^{m(k,T)}a_{s} \le T < \sum_{s=k}^{m(k,T)+1}a_{s}$. From the estimation $\int_{s}^{s+1}\frac{1}{t}\text{d}t<\frac{1}{s}<\int_{s-1}^{s}\frac{1}{t}\text{d}t, s>1$, we have
\begin{align*}
T\ge\sum_{s=k}^{m(k,T)}a_{s}&>\sum_{s=k}^{m(k,T)}\int_{s}^{s+1}\frac{1}{t}\text{d}t\\ &=\int_{k}^{m(k,T)+1}\frac{1}{t}\text{d}t=\ln{\frac{m(k,T)+1}{k}}
\end{align*}
and
\begin{align*}
T<\sum_{s=k}^{m(k,T)+1}a_{s}&<\sum_{s=k}^{m(k,T)+1}\int_{s-1}^{s}\frac{1}{t}\text{d}t \\ &=\int_{k-1}^{m(k,T)+1}\frac{1}{t}\text{d}t=\ln{\frac{m(k,T)+1}{k-1}}.
\end{align*}
From here it follows that $m(k,T)$ satisfies
\begin{align}\label{BoundOfMkT}
(k-1)\exp(T)-1<m(k,T)< k\exp(T)-1.
\end{align}

\begin{lem}\label{PropertyEpsi}
Assume that A1,A2 i), A2 ii), and A3 hold. At the samples $\omega\in\Omega$ where \eqref{SumAkEpijk} holds, for a given $T>0$, for sufficiently large $C>0$ and any $l:l=k,k+1,\cdots,\Big(\big(r(i,m_{i,k}+1)-1\big)\wedge m(k,T)\Big)$, the sequence $\{\epsilon_{i,k}\}$ satisfies
\begin{align}\label{ProEpsi}
 \limsup_{k \to \infty} | \sum_{s=k}^l a_{s}\epsilon_{i,s+1}I_{\{\|\mathbbm{u}_{s}\|\le C\}}|=0\text{ }\forall i \in \mathcal{N},
\end{align}
where $m_{i,k}=\sup\{m:r(i,m) \le k \}$ is the biggest truncation number of agent $i$ up to time $k$.
\end{lem}
\begin{pf}
It is noticed that the inequality $l\le r(i,m_{i,k}+1)-1$ means that at time $l$ the $(m_{i,k}+1)$th truncation has not happened yet for agent $i$, i.e., there is no truncation in \eqref{SigmakHat}-\eqref{AlgorithmSigmak} as $s=k, k+1 \cdots, l$ for agent $i$. By $r(i,0)=1\le k$, the set $\{m:r(i,m) \le k \}$ is nonempty, and then $m_{i,k}$ are well defined. For \eqref{ProEpsi} it suffices to prove
\begin{align}\label{EPOneToFour}
\limsup_{k \to \infty} |\sum_{s=k}^l a_{s}\epsilon_{i,s+1}^{(h)}I_{\{\|\mathbbm{u}_{s}\|\le C\}}|=0,\text{ } h=1,2,3,
\end{align}
where $l=k,k+1,\cdots,\Big(\big(r(i,m_{i,k}+1)-1\big)\wedge m(k,T)\Big)$. From A3, \eqref{EPOneToFour} holds for $h=1$.

By noticing
\begin{align*}
-\epsilon_{i,k+1}^{(2)}=p_{i}\left(y_{i,k+1}-h_{i}(u_{i,k})\right),
\end{align*}
if $l \le k+[\ln k]$, \eqref{EPOneToFour} holds for $h=2$ by \eqref{AkVkToLn} and \eqref{AkYkToLn}. Assume $l > k+[\ln k]$ in what follows. Since the definitions of $h_{i}(\cdot)$ are different for the Hammerstein and Wiener systems, let us consider two cases separately.

1) The dynamics of agent $i$ is a Hammerstein system \eqref{HSystem}.

In this case $\epsilon_{i,k+1}^{(2)}$ can be written as
\begin{align*}
&-\epsilon_{i,k+1}^{(2)}=p_{i}\left(y_{i,k+1}-h_{i}(u_{i,k})\right)\\
&=p_{i}\left(y_{i,k+1}-\frac{d_{i}}{c_{i}}f_{i}(u_{i,k})\right)\\
&=p_{i}\left(y_{i,k+1}-\frac{1}{c_{i}}C_{i}(z)y_{i,k+1}+\frac{1}{c_{i}}D_{i}(z)f_{i}(u_{i,k})-\frac{d_{i}}{c_{i}}f_{i}(u_{i,k})\right) \\
&=\frac{p_{i}}{c_{i}}\Big[ \sum_{v=1}^{p}c_{i,v}(y_{i,k+1}-y_{i,k+1-v}) \\ &+\sum_{v=1}^{q}d_{i,v}\left(f_{i}(u_{i,k-v})-f_{i}(u_{i,k})\right)\Big].
\end{align*}
By taking notice of the discussion for the case of $l\le k+[\ln k]$ given before, for \eqref{EPOneToFour} to hold for $h=2$ it suffices to verify
\begin{align}\label{EPFukSec}
\limsup_{k \to \infty}|\sum_{s=k+[\ln k]+1}^l a_{s}\sum_{v=1}^{q}d_{i,v}(f_{i}(u_{i,s-v})\notag\\-f_{i}(u_{i,s}))I_{\{\|\mathbbm{u}_{s}\|\le C\}}|=0
\end{align}
and
\begin{align}\label{EPYkSec}
\limsup_{k \to \infty}|\sum_{s=k+[\ln k]+1}^l a_{s}\sum_{v=1}^{p}c_{i,v}(y_{i,s+1}\notag\\-y_{i,s+1-v})I_{\{\|\mathbbm{u}_{s}\|\le C\}}|=0
\end{align}
for $l> k+[\ln k]$.

We first prove \eqref{EPFukSec}. Noticing $l\le m(k,T)$ we have
\begin{align}
&|\sum_{s=k+[\ln k]+1}^l a_{s}\sum_{v=1}^{q}d_{i,v}(f_{i}(u_{i,s-v})-f_{i}(u_{i,s}))I_{\{\|\mathbbm{u}_{s}\|\le C\}}| \notag \\
&\le T\sup_{ k+[\ln k]+1\le s \le l}\{\sum_{v=1}^{q}d_{i,v}|f_{i}(u_{i,s-v})-f_{i}(u_{i,s})|I_{\{\|\mathbbm{u}_{s}\|\le C\}}\}.\label{TSupFF}
\end{align}
which means that we need only to prove
\[
\lim_{k\to \infty}|f_{i}(u_{i,s-v})-f_{i}(u_{i,s})|=0,
\]
where $s=k+[\ln k]+1, k+[\ln k]+2,\cdots,l;\text{ } v=1,2,\cdots,q,$ when $|u_{i,s}|\le\|\mathbbm{u}_{s}\|\le C$. By continuity of $f_{i}(\cdot)$ and $|u_{i,s}|\le C$, it is sufficient to show
\begin{align}\label{LimUsvUs}
\lim_{k\to \infty}|u_{i,s}-u_{i,s-v}|=0, \text{ }
\end{align}
where $s=k+[\ln k]+1, k+[\ln k]+2,\cdots,l;\text{ } v=1,2,\cdots,q.$ By noticing $|u_{i,s}-u_{i,s-v}|\le \sum_{h=0}^{v-1}|u_{i,s-h}-u_{i,s-h-1}|$, \eqref{LimUsvUs} follows from
\begin{align}\label{LimUshUsh1}
\lim_{k\to \infty}|u_{i,s-h}-u_{i,s-h-1}|=0,
\end{align}
where $s=k+[\ln k]+1,k+[\ln k]+2,\cdots,l\text{ ; } h=0,1,\cdots,v-1\text{ ; } v=1,2,\cdots,q$. Since $s-h=k+[\ln k]-q+2,k+[\ln k]-q+3,\cdots,l$, for \eqref{LimUshUsh1} it is equivalent to show
\begin{align}\label{LimUh1Uh1}
\lim_{k \to \infty} |u_{i,t+1}-u_{i,t}|=0,
\end{align}
where $ t=k+[\ln k]-q+1,k+[\ln k]-q+2,\cdots,l-1.$

As noticed at the beginning of the proof, there is no truncation in \eqref{SigmakHat}-\eqref{AlgorithmSigmak} as $s=k, \cdots, l$ and $\sigma_{i,k}^{'}=\sigma_{i,k}, u_{i,k}^{'}=u_{i,k}$, i.e.,
\begin{align*}
u_{i,s+1}=u_{i,s}+a_{s}O_{i,s+1}, \text{ }s=k,k+1,\cdots,l-1.
\end{align*}
By \eqref{AkOkToLn} we have
\begin{align}\label{UsPlusOneMinusUs}
|u_{i,s+1}-u_{i,s}|\xrightarrow[k \to \infty]{}0,\text{ } s=k,k+1, \cdots, l-1
\end{align}
which imply \eqref{LimUh1Uh1} and hence \eqref{EPFukSec} holds.

We continue to prove \eqref{EPYkSec}. Similar to \eqref{TSupFF} it is reduced to prove
\begin{align}\label{YS1YS1v}
\lim_{k\to\infty}|y_{i,s+1}-y_{i,s+1-v}|=0,
\end{align}
where $s=k+[\ln k]+1,k+[\ln k]+2,\cdots, l\text{ ; } v=1,2,\cdots,p,$ when $|u_{i,s}|\le\|\mathbbm{u}_{s}\|\le C$. By noticing $(y_{i,s+1}-y_{i,s+1-v})=\sum_{h=0}^{v-1} (y_{i,s-h+1}-y_{i,s-h})$, similar to the proof of \eqref{LimUshUsh1} for \eqref{YS1YS1v} it suffices to prove
\begin{align*}
\lim_{k \to \infty}| y_{i,s-h+1}-y_{i,s-h} |=0
\end{align*}
for $|u_{i,s}|\le\|\mathbbm{u}_{s}\|\le C$, where $s=k+[\ln k]+1,k+[\ln k]+2,\cdots, l\text{ ; } h=0,1,\cdots,v-1\text{ ; } v=1,2,\cdots,p$. Since $s-h=k+[\ln k]-p+2,k+[\ln k]-p+3,\cdots,l$, this is equivalent to verifying
\begin{align}\label{YhPlusOneMinusYh}
\lim_{k \to \infty} | y_{i,k+[\ln{k}]+t+1}-y_{i,k+[\ln k]+t} |=0,
\end{align}
where $t=2-p, 3-p,\cdots,l-k-[\ln k]$, when $|u_{i,s}|\le\|\mathbbm{u}_{s}\|\le C$. From \eqref{EstimateforC}, \eqref{EstForVk} and \eqref{EstForYk}, we have the following chain of inequalities or equalities:
\begin{align}\label{EstYhPlusOneMinusYh}
&|y_{i,k+[\ln{k}]+t+1}-y_{i,k+[\ln{k}]+t}|\notag \\
&\le \| Y_{i,k+[\ln{k}]+t+1}-Y_{i,n_{k}+[\ln{k}]+t} \| \notag \\
&\le \| C_{i}^{[\ln{k}]+t+1}Y_{i,k}-C_{i}^{[\ln{k}]+t}Y_{i,k} \| \notag \\
&+\| \sum_{g=k}^{k+[\ln{k}]+t} C_{i}^{k+[\ln{k}]+t-g}D_{i}f_{i}(u_{i,g})\notag\\&-\sum_{g=k}^{k+[\ln{k}]+t-1} C_{i}^{k+[\ln{k}]+t-1-g}D_{i}f_{i}(u_{i,g})\| \notag \\ & \le re^{-\delta ([\ln{k}]+t)} \alpha_{2}(\ln{k})^{\mu}|re^{-\delta}-1|
+\| C_{i}^{[\ln{k}]+t}D_{i}f_{i}(u_{i,k}) \| \notag \\
&+\|\sum_{g=k}^{k+[\ln{k}]+t-1}C_{i}^{{k}+[\ln{k}]+t-1-g} D_{i}(f_{i}(u_{i,g+1})-f_{i}(u_{i,g})) \| \notag \\
&\le r\alpha_{2}e^{-\delta ([\ln{k}]+t)}(\ln{k})^{\mu} |re^{-\delta}-1|\notag\\
&+r\alpha_{1}\|D_{i}\|e^{-\delta ([\ln{k}]+t)}(\ln k)^{\mu} \notag \\
&+\frac {r\|D_{i}\|}{1-e^{-\delta}}\max_{k \le s \le l-1}| f_{i}(u_{i,s+1})-f(u_{i,s})|.
\end{align}
The first and second terms at the right-hand side of \eqref{EstYhPlusOneMinusYh} tend to zero as $k\to\infty$ and the third term also tends to zero by \eqref{UsPlusOneMinusUs}, $|u_{i,s}|\le\|\mathbbm{u}_{s}\|\le C$, and continuity of $f_{i}(\cdot)$.

2) The dynamics of agent $i$ is a Wiener system \eqref{WSystem}.

In this case $\epsilon_{i,k+1}^{(2)}$ can be written as
\begin{align}
&-\epsilon_{i,k+1}^{(2)}=y_{i,k+1}-h_{i}(u_{i,k})\notag\\
&=f_{i}(v_{i,k+1})-f_{i}(\frac{d_{i}}{c_{i}}u_{i,k})\notag\\
&=f_{i}(v_{i,k+1})-f_{i}\Big(\frac{1}{c_{i}}C_{i}(z)v_{i,k+1}\notag\\ &-\frac{1}{c_{i}}D_{i}(z)u_{i,k}+\frac{d_{i}}{c_{i}}u_{i,k}\Big)\notag\\
&=f_{i}(v_{i,k+1})-f_{i}\Big(v_{i,k+1}+\frac{1}{c_{i}}\sum_{v=1}^{p}c_{i,v}(v_{i,k+1-v}-v_{i,k+1}) \notag\\ &+\frac{1}{c_{i}}\sum_{v=1}^{q}d_{i,v}(u_{i,k}-u_{i,k-v})\Big).
\end{align}
Similar to the preceding proof in 1) it suffices to prove
\begin{align}\label{Ep2Wiener2}
\limsup_{k \to \infty} |\sum_{s=k+[\ln k]+1}^l a_{s}\epsilon_{i,s+1}^{(2)}I_{\{\|\mathbbm{u}_{s}\|\le C\}}|=0
\end{align}
when $l> k+[\ln k]$. From $l\le m(k,T)$ and $\|\mathbbm{u}_{s}\|\le C$ it follows that $\|v_{i,s}\|\le C_{1}$ for some $C_{1}>0$ by \eqref{NorOfVik}. Similar to \eqref{TSupFF}, \eqref{Ep2Wiener2} follows from
\begin{align}\label{LimFVsFVs}
\lim_{k\to\infty}\Big|f_{i}(v_{i,s+1})&-f_{i}\Big(v_{i,s+1}+\frac{1}{c_{i}}\sum_{v=1}^{p} c_{i,v}(v_{i,s+1-v}-v_{i,s+1})\notag \\ &+\frac{1}{c_{i}}\sum_{v=1}^{q}d_{i,v}(u_{i,s}-u_{i,s-v})\Big)\Big|=0,
\end{align}
where $s=k+[\ln k]+1,k+[\ln k]+2,\cdots,l$. By continuity of $f_{i}(\cdot)$, \eqref{LimFVsFVs} follows from
\begin{align}
&\lim_{k\to \infty}|u_{i,s}-u_{i,s-v}|=0, \label{LimUsMUsv}
\end{align}
where $s=k+[\ln k]+1,k+[\ln k]+2,\cdots,l;\text{ } v=1,2,\cdots,q,$ and
\begin{align}
&\lim_{k\to\infty}|v_{i,s+1}-v_{i,s+1-v}|=0, \label{LimVsMVsv}
\end{align}
where $s=k+[\ln k]+1,k+[\ln k]+2,\cdots, l\text{ ; } v=1,2,\cdots,p.$.

Similar to \eqref{LimUsvUs}, we can verify \eqref{LimUsMUsv}, while from \eqref{EstYhPlusOneMinusYh} with $y_{i,k},Y_{i.k},f_{i}(u_{i,k})$ replaced by $v_{i,k},V_{i,k},u_{i,k}$, respectively, we conclude
\begin{align}
\lim_{k\to\infty}|v_{i,k+[\ln k]+t+1}-v_{i,k+[\ln k]+t}|=0, \label{LimVktMVkt}
\end{align}
where $t=2-p,3-p,\cdots,l-k-[\ln k]$. Therefore, \eqref{LimVsMVsv} can be proved in a way similar to that for \eqref{YS1YS1v}.

Similarly, it can be verified that \eqref{EPOneToFour} holds for $h=3$. \qed
\end{pf}

\section{Auxiliary sequences}
Since the growth rate of $f_{i}(\cdot)$ in \eqref{HSystem} and \eqref{WSystem} may be faster than linearly, the SA proposed by \cite{Robbins1951} may diverge. Therefore, we apply the expanding truncation technique introduced in \cite{Chen2002}, and the resulting algorithm is called SAAWET. The distributed algorithm \eqref{SigmakHat} - \eqref{AlgorithmSigmak} is hard to be written as a centralized SAAWET, because the truncation numbers for different agents may be different. By introducing auxiliary sequences, we can transform the DSAAWET to a centralized algorithm by the treatment proposed in \cite{Lei2016}, \cite{Lei2015}, and \cite{Fang2001}.

The discrete-time axis is divided by the sequence $\{r(m):r(m)<\infty\}$ into $[r(m),r(m+1)), m=1,2,\cdots$. Define $\overline{r}(i,m)=r(i,m)\wedge r(m+1)$. When $r(m)<\infty$, by \eqref{RimRmLeDG} we have $0\le\overline{r}(i,m)-r(m)\le d(\mathcal{G})$ which implies that $\overline{r}(i,m)-r(m)$ is bounded. Similar to \cite{Lei2016}, \cite{Lei2015}, and \cite{Fang2001}, we define auxiliary sequences $\{\overline{u}_{i,k}\}$ and $\{\overline{\epsilon}_{i,k+1}\}$ as follows
\begin{align}
\forall k:&r(m) \le k <\overline{r}(i,m),\notag\\
&\begin{cases}\overline{u}_{i,k}=u_{i}^{*},\\ \overline{\epsilon}_{i,k+1}=-\sum_{j \in \mathcal{N}_{i}}p_{ij}(h_{j}(\overline{u}_{j,k})
-h_{i}(u_{i}^{*})),\end{cases} \label{AuxiUBar}\\
\forall k:&\overline{r}(i,m) \le k <r(m+1), \notag\\
&\begin{cases} \overline{u}_{i,k}=u_{i,k},\\
\overline{\epsilon}_{i,k+1}=\epsilon_{i,k+1} +\sum_{j\in\mathcal{N}_{i}}p_{ij}(h_{j}(u_{j,k})-h_{j}(\overline{u}_{j,k})).\end{cases} \label{AuxiSigBar}
\end{align}

For $\{\overline{u}_{i,k}\}$ and $\{\overline{\epsilon}_{i,k+1}\}$ we have the following lemma.
\begin{lem}\label{UBarSigBar}
The sequences $\{\overline{u}_{i,k}\}$ and $\{\overline{\epsilon}_{i,k+1}\}$ satisfy the following recursion:
\begin{align}
\overline{u}_{i,k+1}&=(\overline{u}_{i,k}+a_{k}(g_{i}(\overline{\mathbbm{u}}_{k}) +\overline{\epsilon}_{i,k+1}))\notag\\
&\cdot I_{\{\max_{j\in\mathcal{N}}|\overline{u}_{j,k}+a_{k}(g_{j}(\overline{\mathbbm{u}}_{k}) +\overline{\epsilon}_{j,k+1})| < M_{\overline{\sigma}_{k}}\}} \notag \\
& +u_{i}^{*}I_{\{\max_{j\in\mathcal{N}}|\overline{u}_{j,k}+a_{k}(g_{j}(\overline{\mathbbm{u}}_{k}) +\overline{\epsilon}_{j,k+1})| \ge M_{\overline{\sigma}_{k}}\}},\label{IerationForUBar} \\
\overline{\sigma}_{k+1}&=\overline{\sigma}_{k}+I_{\{\max_{j\in\mathcal{N}}|\overline{u}_{j,k}+a_{k}(g_{j}(\overline{\mathbbm{u}}_{k}) +\overline{\epsilon}_{j,k+1})| \ge M_{\overline{\sigma}_{k}}\}}, \label{IterationForSigBar}
\end{align}
where $\overline{\sigma}_{k}=\max \{\sigma_{j,k},j \in \mathcal{N}\}$.
\end{lem}
\begin{pf}
Define $\tilde{u}_{i,k}, \forall i\in\mathcal{N}$ and $\tilde{\sigma}_{k}$ as follows
\begin{align}
\tilde{u}_{i,k+1}&=(\tilde{u}_{i,k}+a_{k}(g_{i}(\overline{\mathbbm{u}}_{k}) +\overline{\epsilon}_{i,k+1}))\notag\\
&\cdot I_{\{\max_{j\in\mathcal{N}}|\tilde{u}_{j,k}+a_{k}(g_{j}(\overline{\mathbbm{u}}_{k}) +\overline{\epsilon}_{j,k+1})| < M_{\tilde{\sigma}_{k}}\}} \notag \\
& +u_{i}^{*}I_{\{\max_{j\in\mathcal{N}}|\tilde{u}_{j,k} +a_{k}(g_{j}(\overline{\mathbbm{u}}_{k}) +\overline{\epsilon}_{j,k+1})| \ge M_{\tilde{\sigma}_{k}}\}}\label{IerationForUTilde},  \\
\tilde{\sigma}_{k+1}&=\tilde{\sigma}_{k} +I_{\{\max_{j\in\mathcal{N}}|\tilde{u}_{j,k}+a_{k}(g_{j}(\overline{\mathbbm{u}}_{k}) +\overline{\epsilon}_{j,k+1})| \ge M_{\tilde{\sigma}_{k}}\}}, \label{IterationForSigTilde}
\end{align}
where $\tilde{u}_{i,1}=\overline{u}_{i,1}, \tilde{\sigma}_{1}=\overline{\sigma}_{1}$. For \eqref{IerationForUBar} and \eqref{IterationForSigBar}, it suffices to prove $\tilde{u}_{i,k}=\overline{u}_{i,k}, \forall i\in\mathcal{N}, \tilde{\sigma}_{k}=\overline{\sigma}_{k}$, and $\tilde{\sigma}_{k}=\max \{\sigma_{j,k},j \in \mathcal{N}\}\text{ } \forall k\ge1$.

Set $\overline{O}_{i,k+1}=g_{i}(\overline{\mathbbm{u}}_{k}) +\overline{\epsilon}_{i,k+1}$. Consider $k\in[r(m),r(m+1))$. When $k\in[r(m),\overline{r}(i,m))$ from \eqref{AuxiUBar} we have
\begin{align}\label{OBar1}
\overline{O}_{i,k+1}&=\sum_{j\in\mathcal{N}_{i}}p_{ij}(h_{j}(\overline{u}_{j,k}) -h_{i}(\overline{u}_{i,k}))+\overline{\epsilon}_{i,k+1}\notag\\
&=\sum_{j\in\mathcal{N}_{i}}p_{ij}(h_{j}(\overline{u}_{j,k})-h_{i}(u_{i}^{*}))+\overline{\epsilon}_{i,k+1}\notag\\
&=0.
\end{align}
When $k\in[\overline{r}(i,m),r(m+1))$ from \eqref{AuxiSigBar} we have
\begin{align}\label{OBar2}
\overline{O}_{i,k+1}&=\sum_{j\in\mathcal{N}_{i}}p_{ij}(h_{j}(\overline{u}_{j,k})-h_{i}(\overline{u}_{i,k})) +\overline{\epsilon}_{i,k+1}\notag\\
&=\sum_{j\in\mathcal{N}_{i}}p_{ij}(h_{j}(\overline{u}_{j,k})-h_{i}(u_{i,k}))+\overline{\epsilon}_{i,k+1}\notag\\
&=\sum_{j\in\mathcal{N}_{i}}p_{ij}(h_{j}(u_{j,k})-h_{i}(u_{i,k}))+\epsilon_{i,k+1}\notag\\
&=g_{i}(\mathbbm{u}_{k})+\epsilon_{i,k+1}\notag\\
&=O_{i,k+1},
\end{align}
where the last equality is from \eqref{ObservationAgenti}. Combining \eqref{OBar1} and \eqref{OBar2} leads to
\begin{align}\label{OBar}
\overline{O}_{i,k+1}=\begin{cases}0, & k\in[r(m),\overline{r}(i,m)),\\
O_{i,k+1}, &  k\in[\overline{r}(i,m),r(m+1)).\end{cases}
\end{align}

We prove the lemma by induction. Firstly, $\tilde{u}_{i,1}=\overline{u}_{i,1}, \forall i\in\mathcal{N}, \tilde{\sigma}_{1}=\overline{\sigma}_{1}$ and $\tilde{\sigma}_{1}=\max \{\sigma_{j,1},j \in \mathcal{N}\}$ by definitions in \eqref{IerationForUTilde} and \eqref{IterationForSigTilde}. Assume that we have proved $\tilde{u}_{i,s}=\overline{u}_{i,s} \text{ } \forall i\in\mathcal{N}, \tilde{\sigma}_{s}=\overline{\sigma}_{s}$, and $\text{ }\tilde{\sigma}_{s}=\max \{\sigma_{j,s},j \in \mathcal{N}\}$ for $1\le s\le k$. Without loss of generality assume $k\in[r(m),r(m+1))$ for some $r(m)<\infty$. By definition of $r(m)$, there exists $j_{0}$ such that $\sigma_{j_{0},k}=m$ and $\sigma_{j,k}<m+1$ for any $j$. From the inductive assumption it follows that
\begin{align}\label{UTildeK}
\tilde{\sigma}_{k}=\overline{\sigma}_{k} =\max\{\sigma_{j,k},j\in\mathcal{N}\}=m.
\end{align}
Then, for $k+1$ we proceed to prove
\begin{align}
\tilde{u}_{i,k+1}&=\overline{u}_{i,k+1}, \forall i\in\mathcal{N},\label{TiUK1} \\ \tilde{\sigma}_{k+1}&=\overline{\sigma}_{k+1} \text{ ÇÒ } \tilde{\sigma}_{k+1}=\max \{\sigma_{j,k+1},j \in \mathcal{N}\}.\label{TiSigK1}
\end{align}
The proof is carried out by three steps. In the first two steps, we prove \eqref{TiUK1} and \eqref{TiSigK1} for $k\in[r(m),r(m+1)-1)$, respectively. In the third step, we discuss the case where $k=r(m+1)-1$.

Set
\begin{align*}
\mathcal{N}^{(1)}&=\{j\in\mathcal{N}:r(j,m)=r(m)\},\\ \mathcal{N}^{(2)}&=\{j\in\mathcal{N}:r(m)<r(j,m)<r(m+1)\},\\ \mathcal{N}^{(3)}&=\{j\in\mathcal{N}:r(j,m)\ge r(m+1)\}.
\end{align*}
It is clear that $\mathcal{N}=\mathcal{N}^{(1)}\cup\mathcal{N}^{(2)}\cup\mathcal{N}^{(3)}$, and $\mathcal{N}^{(1)}, \mathcal{N}^{(2)}$, and $\mathcal{N}^{(3)}$ are disjoint.

\textbf{Step 1:} We prove that
\begin{align}\label{AbsUtildeLeM}
|\tilde{u}_{j,k}+a_{k}\overline{O}_{j,k+1}|=|\overline{u}_{j,k}+a_{k}\overline{O}_{j,k+1}|
<M_{\tilde{\sigma}_{k}}, \forall j\in\mathcal{N},
\end{align}
i.e., $I_{\{\max_{j\in\mathcal{N}} |\tilde{u}_{j,k}+a_{k}\overline{O}_{j,k+1}| < M_{\tilde{\sigma}_{k}}\}}=1$ for $k\in[r(m),r(m+1)-1)$.

Noticing the inductive assumptions, $\tilde{u}_{j,k}=\overline{u}_{j,k}\text{ } \forall{j}\in\mathcal{N}$, and $\tilde{\sigma}_{k}=\overline{\sigma}_{k}=m$, we verify \eqref{AbsUtildeLeM} for agents in $\mathcal{N}^{(1)}, \mathcal{N}^{(2)}$ and $\mathcal{N}^{(3)}$, respectively.

1) For $j\in\mathcal{N}^{(1)}$, we have $r(j,m)=r(m)$ by definition. Since $k\in[r(m),r(m+1)-1)$, the truncation number of agent $j$ has reached $m$ at time $k$ and the next truncation will not occur at $k+1$. Then it follows that $\sigma_{j,k}=\sigma_{j,r(j,m)}=m=\tilde{\sigma}_{k}$, where the second equality is from the definition of $r(j,m)$, and $\overline{r}(j,m)=r(j,m)\wedge r(m+1)=r(m)$, which means that $[r(m),\overline{r}(j,m))$ in \eqref{AuxiUBar} is empty. By \eqref{SigmakHat} and \eqref{UyiPie}, we have $\sigma_{j,k}^{'}=\sigma_{j,k}$ and $u_{j,k}^{'}=u_{j,k}$, respectively. From \eqref{AuxiSigBar} it follows that $\overline{u}_{j,k}=u_{j,k}$ and from \eqref{OBar} that $\overline{O}_{j,k+1}=O_{j,k+1}$, which implies
\begin{align*}
&|\tilde{u}_{j,k}+a_{k}\overline{O}_{j,k+1}|=|\overline{u}_{j,k}+a_{k}\overline{O}_{j,k+1}|\\
&=|u_{j,k}+a_{k}O_{j,k+1}|<M_{\sigma_{j,k}}=M_{m} =M_{\tilde{\sigma}_{k}},
\end{align*}
where the inequality holds because $k+1<r(m+1)$.

2) For $j\in\mathcal{N}^{(2)}$, we have $r(m)<r(j,m)<r(m+1)$ by definition. Since $k\in[r(m),r(m+1)-1)$, the truncation number of $j$ is smaller than $m$ at $k$. Then it follows that $\sigma_{j,k}\le m=\tilde{\sigma}_{k}$ and $\overline{r}(j,m)=r(j,m)$. By \eqref{AuxiUBar} and \eqref{AuxiSigBar}
\begin{align*}
\overline{u}_{j,k}=\begin{cases}u_{j}^{*},& \text{ } k\in[r(m),r(j,m)), \\
u_{j,k}, & \text{ } k\in[r(j,m),r(m+1)-1).\end{cases}
\end{align*}

From \eqref{OBar} it follows that
\begin{align*}
\overline{O}_{j,k+1}=\begin{cases}0, & k\in[r(m),r(j,m))\\
O_{j,k+1}, & k\in[r(j,m),r(m+1)-1).\end{cases}
\end{align*}
Besides, by \eqref{SigmakHat} and \eqref{UyiPie}, we have $\sigma_{j,k}^{'}=\sigma_{j,k}$ and $u_{j,k}^{'}=u_{j,k}$ for $k\in[r(j,m),r(m+1)-1)$. Therefore,
\begin{align}\label{TilUTwoPart}
&|\tilde{u}_{j,k}+a_{k}\overline{O}_{j,k+1}|=|\overline{u}_{j,k}+a_{k}\overline{O}_{j,k+1}|\notag\\
&=\begin{cases}|u_{j}^{*}|, & k\in[r(m),r(j,m)),\\
|u_{j,k}+a_{k}O_{j,k+1}|, & k\in[r(j,m),r(m+1)-1).\end{cases}
\end{align}
Noticing $|u_{j}^{*}|<M_{0}\le M_{\sigma_{j,k}}$ and $k<r(m+1)-1$, for both cases at the right-hand side of \eqref{TilUTwoPart} we always have
\[
|\tilde{u}_{j,k}+a_{k}\overline{O}_{j,k+1}|<M_{\sigma_{j,k}}\le M_{m}=M_{\tilde{\sigma}_{k}}.
\]

3) For $j\in\mathcal{N}^{(3)}$, we have $r(j,m)>r(m+1)$ by definition. Then we derive $\sigma_{j,k}<m=\tilde{\sigma}_{k}$ and $ \overline{r}(j,m)=r(m+1)$, which means $[\overline{r}(j,m),r(m+1))$ in \eqref{AuxiSigBar} is empty. From \eqref{AuxiUBar} it follows that $\overline{u}_{j,k}=u_{j}^{*}$ and from \eqref{OBar} that $\overline{O}_{j,k+1}=0$. Therefore,
\begin{align*}
|\tilde{u}_{j,k}+a_{k}\overline{O}_{j,k+1}|=|\overline{u}_{j,k}+a_{k}\overline{O}_{j,k+1}|=|u_{j}^{*}|.
\end{align*}
Noticing $|u_{j}^{*}|<M_{0}\le M_{\sigma_{j,k}}$, we have
\[
|\tilde{u}_{j,k}+a_{k}\overline{O}_{j,k+1}|<M_{\sigma_{j,k}}<M_{m}=M_{\tilde{\sigma}_{k}}.
\]
Until now we have proved \eqref{AbsUtildeLeM} for $k\in[r(m),r(m+1)-1)$.

From \eqref{AbsUtildeLeM} it follows that
\[
\tilde{\sigma}_{k+1}=\tilde{\sigma}_{k}=\overline{\sigma}_{k}=\max \{\sigma_{j,k},j \in \mathcal{N}\}
\]
for $k\in[r(m),r(m+1)-1)$. Since $k+1<r(m+1)$, the possible $(m+1)$th truncation can happen only after time $k+1$. Therefore,
\[
\max \{\sigma_{j,k},j \in \mathcal{N}\}=\max \{\sigma_{j,k+1},j \in \mathcal{N}\}=\overline{\sigma}_{i,k+1}
\]
Thus, \eqref{TiSigK1} holds for $k+1$ when $k\in[r(m),r(m+1)-1)$.

\textbf{Step 2:} We proceed to prove that \eqref{TiUK1} holds for $k+1$ when $k\in[r(m),r(m+1)-1)$. Paying attention to
\eqref{IerationForUTilde},\eqref{AbsUtildeLeM} and the inductive assumption, we have
\begin{align}
&\tilde{u}_{i,k+1}=\tilde{u}_{i,k}+a_{k}\overline{O}_{i,k+1},\notag\\
&=\overline{u}_{i,k}+a_{k}\overline{O}_{i,k+1} \text{ , } k\in[r(m),r(m+1)-1)\text{ , }\forall i\in\mathcal{N}.\label{TiUk1TiUk}
\end{align}
Similar to \textbf{Step 1}, we can prove \eqref{TiUK1} for agents in $\mathcal{N}^{(1)}, \mathcal{N}^{(2)}$, and $\mathcal{N}^{(3)}$, respectively.

1) For $i\in\mathcal{N}^{(1)}$, $r(i,m)=r(m)$ by definition, i.e., agent $i$ earlier than other agents reaches $m$ truncations and $\overline{r}(i,m)=r(i,m)=r(m)$. This means that $[r(m),\overline{r}(i,m))$ in
\eqref{AuxiUBar} is empty. By \eqref{SigmakHat}, $\sigma_{i,k}^{'}=\sigma_{i,k}$ and by \eqref{UyiPie}, $u_{i,k}^{'}=u_{i,k}$. From \eqref{AuxiSigBar}, we have $\overline{u}_{i,k}=u_{i,k}$, and also $\overline{u}_{i,k+1}=u_{i,k+1}$ since $k+1<r(m+1)$. From \eqref{OBar} it follows that $\overline{O}_{i,k+1}=O_{i,k+1}$. Since $k+1<r(m+1)$, i.e., there is no truncation for agent $i$ at $k+1$, from \eqref{AlgorithmUk} it follows that
\[
 u_{i,k+1}=u_{i,k}+a_{k}O_{i,k+1},\text{ } k\in[r(m),r(m+1)-1),
\]
which associated with \eqref{TiUk1TiUk} yields $\tilde{u}_{i,k+1}=\overline{u}_{i,k+1}$.

2) For $i\in\mathcal{N}^{(2)}$, $r(m)<r(i,m)<r(m+1)$ by definition, and $\overline{r}(i,m)=r(i,m)$. From \eqref{AuxiUBar} and \eqref{AuxiSigBar} it follows that
\begin{align}\label{Item21}
\overline{u}_{i,k}=\begin{cases}u_{i}^{*}, &\text{ } k\in[r(m),r(i,m)),\\
u_{i,k}, &\text{ } k\in[r(i,m),r(m+1)-1), \end{cases}
\end{align}
and from \eqref{OBar} that
\begin{align}\label{Item22}
\overline{O}_{i,k+1}=\begin{cases}0, &\text{ } k\in[r(m),r(i,m)),\\
O_{i,k+1}, &\text{ } k\in[r(i,m),r(m+1)-1).\end{cases}
\end{align}
Combining \eqref{TiUk1TiUk}, \eqref{Item21} and \eqref{Item22}, we have
\begin{align}\label{Item23}
&\tilde{u}_{i,k+1}\notag\\
&=\begin{cases}u_{i}^{*}, & k\in[r(m),r(i,m)),\\
u_{i,k}+a_{k}O_{i,k+1}=u_{i,k+1}, & k\in[r(i,m),r(m+1)-1).
\end{cases}
\end{align}
We explain the last line of \eqref{Item23}. Notice $k+1<r(m+1)$, i.e., there is no truncation at $k+1$. From \eqref{SigmakHat} and \eqref{UyiPie}, we have $\sigma_{j,k}^{'}=\sigma_{j,k}$ and $u_{j,k}^{'}=u_{j,k}$ for $k\in[r(j,m),r(m+1)-1)$. Then, the last line of \eqref{Item23} is derived from \eqref{AlgorithmUk}.

We now proceed to prove
\begin{align}\label{Item24}
\overline{u}_{i,k+1}=\begin{cases}u_{i}^{*}, & k\in[r(m),r(i,m)),\\
u_{i,k+1}, & k\in[r(i,m),r(m+1)-1).
\end{cases}
\end{align}
It is clear that for $k+1\in[r(m)+1,r(m+1))$ there are three possible cases: $k+1\in(r(m),r(i,m))$, $k+1\in(r(i,m),r(m+1))$, and $k+1=r(i,m)$. For the first two cases \eqref{Item24} is directly derived from
\eqref{AuxiUBar} and \eqref{AuxiSigBar}. We need only to verify the case $k+1=r(i,m)$. To this end, we have $\overline{u}_{i,k+1}=\overline{u}_{i,r(i,m)}=u_{i,r(i,m)}$ by \eqref{AuxiSigBar} and $u_{i,r(i,m)}=u_{i}^{*}$ by the definition of $r(i,m)$. Thus, \eqref{Item24} has been proved. This together with \eqref{Item23} implies that $\tilde{u}_{i,k+1}=\overline{u}_{i,k+1}$ for $k\in[r(m),r(m+1)-1)$.

3) For $i\in\mathcal{N}^{(3)}$, $r(i,m)\ge r(m+1)$ by definition and hence $\overline{r}(i,m)=r(m+1)$. This means that $[\overline{r}(j,m),r(m+1))$ in \eqref{AuxiSigBar} is empty. By \eqref{AuxiUBar}, $\overline{u}_{i,k}=u_{i}^{*}$ and $\overline{u}_{i,k+1}=u_{i}^{*}$ for $k+1\in(r(m),r(m+1))$. From \eqref{OBar}, we have $\overline{O}_{i,k+1}=0$ which associated with \eqref{TiUk1TiUk} leads to $\tilde{u}_{i,k+1}=u_{i}^{*}=\overline{u}_{i,k+1}$.

\textbf{Step 3:} If $r(m+1)=\infty$, then we have completed the proof. This is because it is unnecessary to consider the case $k=r(m+1)-1$ . It remains to show \eqref{TiUK1} and \eqref{TiSigK1} when $r(m+1)<\infty$ and $k=r(m+1)-1$.

For $k=r(m+1)-1$, at $k+1$ the $(m+1)$th truncation occurs for some agent $j_{0}$, i.e., there exists $j_{0}$ such that $|u_{j_{0},k}+a_{k}O_{j_{0},k+1}|\ge M_{m}$ and then $\sigma_{j_{0},k+1}=m+1$. Since $r(j_{0},m)<k+1=r(m+1)$, we have $\overline{r}(j_{0},m)=r(j_{0},m)$. From \eqref{AuxiSigBar} it follows that $\overline{u}_{j_{0},k}=u_{j_{0},k}$ and from \eqref{OBar} that $\overline{O}_{j_{0},k+1}=O_{j_{0},k+1}$. Therefore, by the inductive assumption we have
\begin{align}
\tilde{\sigma}_{k+1}&=\tilde{\sigma}_{k} +I_{\{\max_{j\in\mathcal{N}}|\tilde{u}_{j,k}+a_{k}\overline{O}_{j,k+1}| \ge M_{\tilde{\sigma}_{k}}\}}\notag\\
&=\tilde{\sigma}_{k} +I_{\{|u_{j_{0},k}+a_{k}O_{j_{0},k+1}| \ge M_{\tilde{\sigma}_{k}}\}}\notag\\
&=\overline{\sigma}_{k}+1.\label{TilSigK1}
\end{align}
Further by \eqref{UTildeK}, the chain of equalities \eqref{TilSigK1} can be continued as follows
\begin{align*}
&\tilde{\sigma}_{k+1}=\overline{\sigma}_{k}+1=m+1=\sigma_{j_{0},k+1}\\
&=\overline{\sigma}_{r(m+1)}
=\max \{\sigma_{j,k+1},j \in \mathcal{N}\}=\overline{\sigma}_{i,k+1}.
\end{align*}
Thus, \eqref{TiSigK1} is proved for $k=r(m+1)-1$. We then have $\tilde{u}_{i,k+1}=u_{i}^{*}\text{ } \forall i\in\mathcal{N}$ from \eqref{IerationForUTilde}, and $\overline{u}_{i,k+1}=\overline{u}_{i,r(m+1)}=u_{i}^{*}\text{ }\forall i\in\mathcal{N}$ from \eqref{AuxiUBar}. This means that $\tilde{u}_{i,k+1}=\overline{u}_{i,k+1} \text{ }\forall i\in\mathcal{N}, k=r(m+1)-1$.

Thus, we have proved \eqref{TiUK1} and \eqref{TiSigK1} for $k+1$ and completed the proof. \qed
\end{pf}

Set $\overline{\mathbbm{u}}_{k}=[\overline{u}_{1,k},\cdots,\overline{u}_{N,k}]^T, \overline{\epsilon}_{k+1}=[\overline{\epsilon}_{1,k+1},\cdots, \overline{\epsilon}_{N,k+1}]^T$ and $ \mathbbm{g}(\overline{\mathbbm{u}})=[g_{1}(\overline{\mathbbm{u}}),\cdots, g_{N}(\overline{\mathbbm{u}})]^{T}$. Then \eqref{IerationForUBar} and \eqref{IterationForSigBar} can be written as
\begin{align}
\overline{\mathbbm{u}}_{k+1}&=(\overline{\mathbbm{u}}_{k} +a_{k}(\mathbbm{g}(\overline{\mathbbm{u}}_{k}) +\overline{\epsilon}_{k+1}))\notag\\
&\cdot I_{\{\|\overline{\mathbbm{u}}_{k} +a_{k}(\mathbbm{g}(\overline{\mathbbm{u}}_{k})+\overline{\epsilon}_{k+1}) \|_{\infty} < M_{\overline{\sigma}_{k}}\}} \notag \\
&+u^{*}I_{\{\|\overline{\mathbbm{u}}_{k}+a_{k}(\mathbbm{g}(\overline{\mathbbm{u}}_{k}) +\overline{\epsilon}_{k+1}) \|_{\infty} \ge M_{\overline{\sigma}_{k}}\}},\label{CompatForUBar}\\
\overline{\sigma}_{k+1}&=\overline{\sigma}_{k}+I_{\{\|\overline{\mathbbm{u}}_{k} +a_{k}(\mathbbm{g}(\overline{\mathbbm{u}}_{k})+\overline{\epsilon}_{k+1}) \|_{\infty} \ge M_{\overline{\sigma}_{k}}\}}.\label{CompatForSigmaBar}
\end{align}

\begin{lem}\label{PropertyEpsiBar}
Assume A1, A2 i), A2 ii), A3, and A4 hold. At the samples $\omega \in \Omega$ where \eqref{SumAkEpijk} holds, for a given $T>0$ and for sufficiently large $C>0$, the noise $\{\overline{\epsilon}_{i,k+1}\}$ defined in \eqref{AuxiUBar} and \eqref{AuxiSigBar} satisfies
\begin{align}\label{AsEpBars}
\lim_{k \to \infty}|&\sum_{s=k}^{l} a_{s}\overline{\epsilon}_{i,s+1}I_{\{\|\overline{\mathbbm{u}}_{s}\|\le C\}}|=0,\text{ }\forall i \in \mathcal{N}, \notag\\ &l=k,k+1,\cdots,\Big(\big(r(m_{k}+1)-1\big)\wedge m(k,T)\Big),
\end{align}
where $m_{k}=\sup\{m:r(m) \le k\}$ is the biggest number of truncations occurred at agents up to time $k$.
\end{lem}
\begin{pf}
It is worth noting that the inequality or equality $l\le r(m_{k}+1)-1$ means that the $(m_{k}+1)$th truncation has not happened yet for any agent, i.e., there is no truncation in \eqref{SigmakHat}-\eqref{AlgorithmSigmak} as $s=k,k+1,\cdots,l$ for any agent.

Similar to Lemma~\ref{PropertyEpsi}, since $r(0)=1\le k$, the set $\{m:r(m) \le k\}$ is nonempty. This means that $m_{k}$ is well defined. Besides, we have
\begin{align}\label{SumEpsiBar}
&|\sum_{s=k}^{l} a_{s}\overline{\epsilon}_{i,s+1}I_{\{\|\overline{\mathbbm{u}}_{s}\|\le C\}}|\le |\sum_{s=k}^{l} a_{s}\overline{\epsilon}_{i,s+1}\notag\\
&\cdot I_{\{r(m_{k})\le s<\overline{r}(i,m_{k})\}}I_{\{\|\overline{\mathbbm{u}}_{s}\|\le C\}}|\notag\\
&+|\sum_{s=k}^{l} a_{s}\overline{\epsilon}_{i,s+1}I_{\{\overline{r}(i,m_{k})\le s < r(m_{k}+1)\}}I_{\{\|\overline{\mathbbm{u}}_{s}\|\le C\}}|.
\end{align}

First, let us consider the first term at the right-hand side of \eqref{SumEpsiBar}. By noticing \eqref{AuxiUBar} it follows that
\begin{align}\label{SumEpsiBar1}
&|\sum_{s=k}^{l} a_{s}\overline{\epsilon}_{i,s+1}I_{\{r(m_{k})\le s<\overline{r}(i,m_{k})\}}I_{\{\|\overline{\mathbbm{u}}_{s}\|\le C\}}|\notag\\
&=|\sum_{s=k}^{l}a_{s}\sum_{j \in \mathcal{N}_{i}}p_{ij}(h_{j}(\overline{u}_{j,k})
-h_{i}(u_{i}^{*}))\notag\\
\cdot &I_{\{r(m_{k})\le s<\overline{r}(i,m_{k})\}}I_{\{\|\mathbbm{u}_{s}\|\le C\}}|.
\end{align}
From \eqref{AuxiUBar} and \eqref{AuxiSigBar} we have
\begin{align}\label{UBarjk}
\overline{u}_{j,k}=\begin{cases}u_{j}^{*}, &k\in[r(m),\overline{r}(j,m)),\\
u_{j,k}, & k\in[\overline{r}(j,m),r(m+1)).
\end{cases}
\end{align}
Therefore, similar to \eqref{EstUik}, we see $|\overline{u}_{j,k}|\le\ln(k-1+c_{M})$ for sufficiently large $k$. This together with the boundedness of $\overline{r}(i,m_{k})-r(m_{k})$ as proved for \eqref{AkGkToLn} implies that the right-hand side of \eqref{SumEpsiBar1} tends to zero as $k\to\infty$.

We now analyze the second term at the right-hand side of \eqref{SumEpsiBar}. From \eqref{AuxiSigBar} it follows that
\begin{align}\label{SumEpsiBar2}
&|\sum_{s=k}^{l} a_{s}\overline{\epsilon}_{i,s+1}I_{\{\overline{r}(i,m_{k})\le  s < r(m_{k}+1)\}}I_{\{\|\overline{\mathbbm{u}}_{s}\|\le C\}}|\notag\\
&\le |\sum_{s=k}^{l}a_{s}\epsilon_{i,s+1}I_{\{\overline{r}(i,m_{k})\le s < r(m_{k}+1)\}}I_{\{\|\mathbbm{u}_{s}\|\le C\}}|\notag\\
&+|\sum_{s=k}^{l}a_{s}\sum_{j\in\mathcal{N}_{i}}p_{ij}(h_{j}(u_{j,s})-h_{j}(\overline{u}_{j,s}))\notag\\
&\cdot I_{\{\overline{r}(i,m_{k})\le s < r(m_{k}+1)\}}I_{\{\|\mathbbm{u}_{s}\|\le C\}}|.
\end{align}
When $\overline{r}(i,m_{k})=r(i,m_{k})$, by Lemma~\ref{PropertyEpsi} the first term at the right-hand side of \eqref{SumEpsiBar2} tends to zero as $k\to\infty$, while it is zero when $\overline{r}(i,m_{k})=r(m_{k}+1)$.

In view of \eqref{UBarjk}, we have
\begin{align*}
&\sum_{s=k}^{l}a_{s}\sum_{j\in\mathcal{N}_{i}}p_{ij}(h_{j}(u_{j,s})-h_{j}(\overline{u}_{j,s})) \notag\\ &=\sum_{j\in\mathcal{N}_{i}}p_{ij}\sum_{s=k}^{l}a_{s}(h_{j}(u_{j,s})-h_{j}(\overline{u}_{j,s}))\notag\\
&=\sum_{j\in\mathcal{N}_{i}}p_{ij}\sum_{s=k}^{l}a_{s}(h_{j}(u_{j,s})-h_{j}(\overline{u}_{j,s}))I_{\{r(m_{k})\le s < \overline{r}(j,m_{k})\}}\notag\\
&=\sum_{j\in\mathcal{N}_{i}}p_{ij}\sum_{s=k}^{l}a_{s}(h_{j}(u_{j,s})-h_{j}(u_{j}^{*}))I_{\{r(m_{k})\le s < \overline{r}(j,m_{k})\}}.
\end{align*}
Therefore, the second term at the right-hand side of \eqref{SumEpsiBar2} can be rewritten as
\begin{align}\label{SumEpsiBar3}
&|\sum_{s=k}^{l}a_{s}\sum_{j\in\mathcal{N}_{i}}p_{ij}(h_{j}(u_{j,s})-h_{j}(\overline{u}_{j,s}))\notag\\ &\cdot I_{\{\overline{r}(i,m_{k})\le s < r(m_{k}+1)\}}I_{\{\|\mathbbm{u}_{s}\|\le C\}}|\notag\\
&=|\sum_{j\in\mathcal{N}_{i}}p_{ij}\sum_{s=k}^{l}a_{s}(h_{j}(u_{j,s})-h_{j}(u_{j}^{*})) \notag\\
&\cdot I_{\{\overline{r}(i,m_{k})\le s < r(m_{k}+1)\}}I_{\{r(m_{k})\le s < \overline{r}(j,m_{k})\}}I_{\{\|\mathbbm{u}_{s}\|\le C\}}|.
\end{align}
By the boundedness of $\overline{r}(j,m_{k})-r(m_{k})$, similar to \eqref{AkGkToLn}, the right-hand side of \eqref{SumEpsiBar3} tends zero as $k\to\infty$.

Combining \eqref{SumEpsiBar}, \eqref{SumEpsiBar1}, \eqref{SumEpsiBar2}, and \eqref{SumEpsiBar3} we conclude the lemma. \qed
\end{pf}

The sequence $\{\overline{\mathbbm{u}}_{k}\}$ has the following property along any its convergent subsequence.
\begin{lem}\label{PropertyConSubSeq}
Assume A1, A2 i), A2 ii), A3, and A4 hold. Let $\{\overline{\mathbbm{u}}_{k}\}$ be generated by \eqref{CompatForUBar} and \eqref{CompatForSigmaBar} and let $\{\overline{\mathbbm{u}}_{n_{k}}\} \subset \{\overline{\mathbbm{u}}_{k}\}$ be a its convergent subsequence. For any sample $\omega \in \Omega$ where \eqref{SumAkEpijk} holds, there exist $c_{1}>0$ and $T_{1}>0$ such that
\begin{align}\label{UsPlusoneMinusUnk}
\| \overline{\mathbbm{u}}_{s+1}-\overline{\mathbbm{u}}_{n_{k}} \| \le c_{1}T,\forall s:n_{k}\le s \le m(n_{k},T) \text{ }\forall T \in [0,T_{1}]
\end{align}
for sufficiently large $k$.
\end{lem}
\begin{pf}
By definition $\mathbbm{g}(\cdot):\mathbb{R}^N \to \mathbb{R}^N$ is continuous. By the boundedness of $\{\overline{\mathbbm{u}}_{n_{k}}\}$, take $M>\sup_{k}\{\|\mathbbm{u}_{n_{k}}\|\}$ and set $G=\sup_{\|\mathbbm{u}\| \le M+1}\|\mathbbm{g}(\mathbbm{u})\|$. Take $c_{1}=2NG$ and $T_{1}>0$ such that $c_{1}T_{1} \le 1$. We now show that $c_{1}$ and $T_{1}$ satisfy \eqref{UsPlusoneMinusUnk} for sufficiently large $k$.

For $T \in [0,T_{1}]$, set
\begin{align}\label{DefForSk}
s_{k}\triangleq\min \{s \ge n_{k},\|\overline{\mathbbm{u}}_{s+1}-\overline{\mathbbm{u}}_{n_{k}}\| > c_{1}T \}.
\end{align}
We need only to prove $s_{k} > m(n_{k},T)$ for sufficiently large $k$. If $s_{k}=\infty$, i.e., $\{s \ge n_{k},\|\overline{\mathbbm{u}}_{s+1}-\overline{\mathbbm{u}}_{n_{k}}\| > c_{1}T \}$ is empty, then we have completed the proof. In the following we assume $s_{k}<\infty$.

We first prove
\begin{align}\label{SkPlusOneLe}
s_{k} \le r(\overline{\sigma}_{n_{k}}+1)-1
\end{align}
for sufficiently large $k$. If $\lim_{k \to \infty}\overline{\sigma}_{k}=\sigma<\infty$, i.e., the largest number of truncations the algorithm can reach is $\sigma$, then $r(\overline{\sigma}_{n_{k}}+1)=r(\sigma+1)=\inf\emptyset=\infty$ for sufficiently large $k$. This implies \eqref{SkPlusOneLe}. On the other hand, if $\lim_{k \to \infty}\overline{\sigma}_{k}=\infty$, then $M_{\overline{\sigma}_{n_{k}}}>M+1$ for sufficiently large $k$. Since
\[
\|\overline{\mathbbm{u}}_{s+1}-\overline{\mathbbm{u}}_{n_{k}}\| \le c_{1}T, \text{ }\forall s:n_{k} \le s < s_{k},
\]
we have
\begin{align}\label{EstForUBarsPlusOne}
\|\overline{\mathbbm{u}}_{s+1}\| \le \|\overline{\mathbbm{u}}_{n_{k}}\|+c_{1}T < M+1,\forall s:n_{k} \le s < s_{k}.
\end{align}
This means that there is no truncation in \eqref{CompatForUBar} and \eqref{CompatForSigmaBar} from $n_{k}$ to $s_{k}$ and thus \eqref{SkPlusOneLe} is proved.

We now show
\begin{align}\label{SkBiggerMkT}
s_{k} > m(n_{k},T).
\end{align}
Assume the converse $s_{k} \le m(n_{k},T)$. By \eqref{EstForUBarsPlusOne} and the definition of $G$, we have
\begin{align}\label{AsGBars}
|\sum_{s=n_{k}}^{s_{k}}a_{s}\overline{g}_{i}(\overline{\mathbbm{u}}_{s})| \le GT=\frac{c_{1}T}{2N}.
\end{align}
Moreover, from \eqref{SkPlusOneLe} and Lemma~\ref{PropertyEpsiBar} (where $C>M+1$), it follows that
\begin{align}\label{AsEpBarsPlusOne}
|\sum_{s=n_{k}}^{s_{k}}a_{s}\overline{\epsilon}_{i,s+1}| \le \frac{c_{1}T}{2N}
\end{align}
for sufficiently large $k$. Thus, combining \eqref{AsGBars} and \eqref{AsEpBarsPlusOne} we have
\begin{align*}
&\|\overline{\mathbbm{u}}_{s_{k}+1}-\overline{\mathbbm{u}}_{n_{k}}\| \le \sum_{i=1}^N |\overline{u}_{i,s_{k}+1}-\overline{u}_{i,n_{k}}|  \notag \\
&\le \sum_{i=1}^N (|\sum_{s=n_{k}}^{s_{k}}a_{s}\overline{g}_{i}(\overline{\mathbbm{u}}_{s})|+
|\sum_{s=n_{k}}^{s_{k}}a_{s}\overline{\epsilon}_{i,s+1}|) \le c_{1}T,
\end{align*}
which contradicts with \eqref{DefForSk}. So, the converse assumption is false, and $s_{k} > m(n_{k},T)$ for sufficiently large $k$. Thus, the lemma is proved. \qed
\end{pf}

\section{Convergence of algorithm}
Noticing
\begin{align*}
\mathbbm{g}(\overline{\mathbbm{u}}_{k})=(P-D)\mathbbm{h}(\overline{\mathbbm{u}}_{k}) =-L\mathbbm{h}(\overline{\mathbbm{u}}_{k}),
\end{align*}
we can rewrite \eqref{CompatForUBar} and \eqref{CompatForSigmaBar} as
\begin{align}
\overline{\mathbbm{u}}_{k+1}&=\left(\overline{\mathbbm{u}}_{k}+a_{k}(-L\mathbbm{h}(\overline{\mathbbm{u}}_{k}) +\overline{\epsilon}_{k+1})\right)\notag\\
&\cdot I_{\{ \|\overline{\mathbbm{u}}_{k}+a_{k}(-L\mathbbm{h}(\overline{\mathbbm{u}}_{k})+\overline{\epsilon}_{k+1}) \|_{\infty} \le M_{\overline{\sigma}_{k}}\}} \notag \\ &+\mathbbm{u}^{*}I_{\{\|\overline{\mathbbm{u}}_{k}+a_{k}(-L\mathbbm{h}(\overline{\mathbbm{u}}_{k})+\overline{\epsilon}_{k+1}) \|_{\infty} > M_{\overline{\sigma}_{k}}\}},\label{UkBarHUkBar}\\
\overline{\sigma}_{k+1}&=\overline{\sigma}_{k}
+I_{\{\|\overline{\mathbbm{u}}_{k}+a_{k}(-L\mathbbm{h}(\overline{\mathbbm{u}}_{k})+\overline{\epsilon}_{k+1}) \|_{\infty} > M_{\overline{\sigma}_{k}}\}}.\label{SigmaBarHUkBar}
\end{align}

\begin{lem}\label{SingleRoot}
Assume A2 iii) holds. The set $J\cap\{\mathbbm{u}\in\mathbb{R}^{N}: \mathbbm{1}^{T}\mathbbm{u}=c\}$ is a singleton, where $c$ is a constant.
\end{lem}
\begin{pf}
In fact, under condition A2 iii) for any given $\alpha_{i}>0,i=1,2,\cdots,N$ and constant $c$, the system of equations with respect to $u_{1},\cdots,u_{N}$
\begin{align}\label{PropertyMotonFuc}
\begin{cases}
h_{1}(u_{1})=\cdots=h_{N}(u_{N}),\\
\alpha_{1} u_{1}+\cdots+\alpha_{N} u_{N}=c
\end{cases}
\end{align}
has a unique solution. Set $h_{1}(u_{1})=\cdots=h_{N}(u_{N})=b$. By A2 iii) $h_{i}(\cdot)$ has the inverse function $h_{i}^{-1}(\cdot)$, which is strictly monotonically increasing with range $(-\infty,+\infty)$. For any given constant $c$, consider the equation of $b$
\begin{align}
\alpha_{1} h_{1}^{-1}(b)+\cdots+\alpha_{N} h_{N}^{-1}(b)=c. \notag
\end{align}
Since $\alpha_{i} > 0$, $\alpha_{1} h_{1}^{-1}(\cdot)+\cdots+\alpha_{N} h_{N}^{-1}(\cdot)$ is also strictly monotonically increasing with range $(-\infty,+\infty)$. Therefore, the above equation has a unique solution denoted as $b_{c}$, which means the original system of equations \eqref{PropertyMotonFuc} has the unique solution $u_{i}=h_{i}^{-1}(b_{c}),i=1,2,\cdots,N$. Thus the lemma is proved, since $\mathbbm{1}^{T}\mathbbm{u}=u_{1}+\cdots+u_{N}$ and the first $N-1$ equations in \eqref{PropertyMotonFuc} are equivalent to $h(u)\in\text{Span}\{\mathbbm{1}\}$. \qed
\end{pf}

Notice that A2 i) and iii) imply
\begin{align}\label{PropertyForIntHi}
0 \le \int_{u_{i}^{(0)}}^{c} h_{i}(t) \text{d}t \xrightarrow[|c| \to \infty]{} \infty,
\end{align}
where $u_{i}^{(0)}$ is the root of $h_{i}(\cdot)$. Define the continuously differentiable function $v(\cdot):\mathbb{R}^N \to \mathbb{R}$ as
\begin{align}\label{DefForVu}
v(\mathbbm{u})=\sum_{i \in \mathcal{N}} \int_{u_{i}^{(0)}}^{u_{i}}h_{i}(t) \text{d}t.
\end{align}
It is clear that $\nabla_{\mathbbm{u}}v(\mathbbm{u})=\mathbbm{h}(\mathbbm{u})$. In view of \eqref{PropertyForIntHi} and \eqref{DefForVu}, there exists a constant $c_{0}>0$ such that $c_{0} > \|\mathbbm{u}^*\|_{\infty}$ and
\begin{align}\label{ChoiseForC0}
& \min \{ \int_{u_{i}^{(0)}}^{c_{0}}h_{i}(t) \text{d}t , \int_{u_{i}^{(0)}}^{-c_{0}}h_{i}(t) \text{d}t\} \notag \\
& > \sum_{j \in \mathcal{N}}\int_{u_{j}^{(0)}}^{u_{j}^*} h_{j}(t)\text{d}t=v(\mathbbm{u}^{*}).
\end{align}
By A4, the Laplace matrix $L$ of the undirected graph has the following property
\begin{align*}
\mathbbm{y}^TL\mathbbm{y}=\frac{1}{2}\sum_{(i,j) \in \mathcal{E}}p_{ij}(y_{i}-y_{j})^2=\frac{1}{2}\sum_{i,j \in \mathcal{N}}p_{ij}(y_{i}-y_{j})^2,
\end{align*}
for $\mathbbm{y}=[y_{1},\cdots,y_{N}]^{T} \in \mathbb{R}^{N}$. Therefore, for any $0<\delta<\Delta$ there exists $\zeta > 0$ such that
\begin{align}\label{NablaVPIHUZeta}
&\sup_{\delta \le d(\mathbbm{u},J) \le \Delta}\nabla_{\mathbbm{u}}^T v(\mathbbm{u})(-L)\mathbbm{h}(\mathbbm{u})\notag\\
&=\frac{1}{2}\sup_{\delta \le d(\mathbbm{u},J) \le \Delta} -\sum_{(i,j) \in \mathcal{E}}p_{i,j}(h_{i}(u_{i})-h_{j}(u_{j}))^2 \le -\zeta.
\end{align}

Theorem~2.2.1 in \cite{Chen2002} is not applicable for proving the convergence $d(\overline{\mathbbm{u}}_{k},J)\xrightarrow[k\to\infty]{}0$ for the algorithm \eqref{UkBarHUkBar} and \eqref{SigmaBarHUkBar}, since we cannot guarantee that $v(J)\triangleq\{v(\mathbbm{u})|\mathbbm{u}\in J\}$ is nowhere dense. However, the idea of proof of that theorem can still be used.

Define $S^*(\epsilon)\triangleq\{\mathbbm{u}: \mathbbm{1}_{N}^T\mathbbm{u} \in [\mathbbm{1}_{N}^T\mathbbm{u}^*-\epsilon,\mathbbm{1}_{N}^T\mathbbm{u}^*+\epsilon]\},\epsilon \ge 0$ and $S^*(0)\triangleq\{u: \mathbbm{1}_{N}^T\mathbbm{u}=\mathbbm{1}_{N}^T\mathbbm{u}^*\}$. To prove the convergence of SAAWET the key step is to show that the truncation ceases in a finite number of steps. This also is the key issue for convergence analysis of the algorithm \eqref{UkBarHUkBar} and \eqref{SigmaBarHUkBar}, and its proof is motivated by that for SAAWET given in \cite{Chen2002}.

\begin{lem}\label{FiniteTruncation}
Assume A1-A4 hold. For samples $\omega \in \Omega$ where \eqref{SumAkEpijk} holds,
\begin{align}\label{FinitenessForSigmakBar}
\lim_{k \to \infty}\overline{\sigma}_{k} \triangleq \overline{\sigma} < \infty.
\end{align}
\end{lem}
\begin{pf}
Assume the converse $\overline{\sigma}_{k} \xrightarrow[k \to \infty]{} \infty$, i.e., there exists a sequence $\{n_{l}\}_{l\ge1}$ such that $\overline{\sigma}_{n_{l}}=\overline{\sigma}_{n_{l}-1}+1$ and $\overline{\mathbbm{u}}_{n_{l}}=\mathbbm{u}^*$. To reach a contradiction we complete the proof by three steps.

\textbf{Step 1:} We first show
\begin{align}\label{DistanceForUkBarSstar0}
d(\overline{\mathbbm{u}}_{k},S^*(0)) \xrightarrow[k \to \infty]{} 0
\end{align}
under the converse assumption.

The discrete-time axis can be divided into $[n_{l},n_{l+1}), l=1,2,\cdots$ by $\{n_{l}\}$ and we investigate the limit of $\overline{\mathbbm{u}}_{k}$ on these intervals. Since $\overline{\mathbbm{u}}_{n_{l}}=\mathbbm{u}^*$, \eqref{DistanceForUkBarSstar0} is satisfied at $\overline{\mathbbm{u}}_{n_{l}}$. We need only to consider $(n_{l},n_{l+1})$. When $T>0$ is sufficiently small, by Lemma~\ref{PropertyConSubSeq} there is no truncation for \eqref{UkBarHUkBar} and \eqref{SigmaBarHUkBar} from $n_{l}$ to $m(n_{l},T)$ for sufficiently large $l$, i.e., $m(n_{l},T)< r(\overline{\sigma}_{n_{l}}+1)=n_{l+1}$. This means $r(\overline{\sigma}_{n_{l}})=n_{l}\le m(n_{l},T)<r(\overline{\sigma}_{n_{l}}+1)=n_{l+1}$. From \eqref{BoundOfMkT} it follows that $m(n_{l},T)-n_{l}>(\exp(T)-1)n_{l}-\exp(T)-1\xrightarrow[l\to\infty]{}\infty$. Let $i_{l}\in \text{arg}\max_{j\in \mathcal{N}}\{r(j,\overline{\sigma}_{n_{l}})\}$. By $r(i_{l},\overline{\sigma}_{n_{l}})-r(\overline{\sigma}_{n_{l}})<d(\mathcal{G})$, we have
\[
n_{l}\le r(i_{l},\overline{\sigma}_{n_{l}}) \le m(n_{l},T)<r(\overline{\sigma}_{n_{l}}+1)=n_{l+1}
\]
for sufficiently large $l$. In the following, the limit of $\overline{\mathbbm{u}}_{k}$ will be shown on $(n_{l},r(i_{l},\overline{\sigma}_{n_{l}})]$ and $(r(i_{l},\overline{\sigma}_{n_{l}}),n_{l+1})$, respectively.

1) For $k\in(n_{l}, r(i_{l},\overline{\sigma}_{n_{l}})]$
\begin{align*}
\overline{\mathbbm{u}}_{k}=\overline{\mathbbm{u}}_{k-1} +a_{k-1}(-L\mathbbm{h}(\overline{\mathbbm{u}}_{k-1})
+\overline{\epsilon}_{k}).
\end{align*}
Multiplying this equation with $\mathbbm{1}_{N}^T$ from the left and noticing that $\mathbbm{1}_{N}$ is the left eigenvector of $L$ corresponding to the eigenvalue $0$ and that $\overline{\mathbbm{u}}_{n_{l}}=\mathbbm{u}^*$, we have
\begin{align}
\mathbbm{1}_{N}^T\overline{\mathbbm{u}}_{k}&=\mathbbm{1}_{N}^T\overline{\mathbbm{u}}_{k-1} +a_{k-1}\mathbbm{1}_{N}^T\overline{\epsilon}_{k} \notag \\
&=\mathbbm{1}_{N}^T \mathbbm{u}^*+\sum_{s=n_{l}}^{k-1}a_{s}\mathbbm{1}_{N}^T\overline{\epsilon}_{s+1}. \label{MulLef1}
\end{align}
Since $k-1<m(n_{l},T)$, by \eqref{EstForUBarsPlusOne}, and Lemma~\ref{PropertyEpsiBar} (where $C>M+1$) the second term at the right-hand side of \eqref{MulLef1} tends to zero as $n_{l} \to \infty$. Thus, $\mathbbm{1}_{N}^{T}\overline{\mathbbm{u}}_{k}\xrightarrow[k\to\infty]{} \mathbbm{1}_{N}^{T}\mathbbm{u}^{*}$ for $k\in (n_{l},r(i_{l},\overline{\sigma}_{n_{l}})]$.

2) For $k\in(r(i_{l},\overline{\sigma}_{n_{l}}),n_{l+1})$, noticing $ r(j,\overline{\sigma}_{n_{l}})< k <r(\overline{\sigma}_{n_{l}}+1), \forall j$, by \eqref{ObservationAgenti} and
\eqref{OBar}, it follows that $\overline{\mathbbm{O}}_{k}=\mathbbm{O}_{k}=(P-D)\mathbbm{y}_{k} +\epsilon_{k}^{(1)}$, where $\mathbbm{O}_{k}=[O_{1,k},\cdots,O_{N,k}]^{T}$ and $\epsilon_{k}^{(1)}=[\epsilon_{1,k}^{(1)},\cdots,\epsilon_{N,k}^{(1)}]$. Then, we have
\begin{align*}
&\overline{\mathbbm{u}}_{k}=\overline{\mathbbm{u}}_{k-1}+a_{k-1}\overline{\mathbbm{O}}_{k}\\ &=\overline{\mathbbm{u}}_{k-1}+a_{k-1}\mathbbm{O}_{k}=\overline{\mathbbm{u}}_{k-1} +a_{k-1}(-L\mathbbm{y}_{k} +\epsilon_{k}^{(1)}).
\end{align*}
Multiplying this equation with $\mathbbm{1}_{N}^T$ from left and noticing $\mathbbm{1}_{N}$ is the left eigenvector of $L$ corresponding to eigenvalue $0$, we derive
\begin{align}
\mathbbm{1}_{N}^T\overline{\mathbbm{u}}_{k}&=\mathbbm{1}_{N}^T\overline{\mathbbm{u}}_{k-1} +a_{k-1}\mathbbm{1}_{N}^T\epsilon_{k}^{(1)} \notag \\
&=\mathbbm{1}_{N}^T \overline{\mathbbm{u}}_{r(i_{l},\overline{\sigma}_{n_{l}})} +\sum_{s=r(i_{l},\overline{\sigma}_{n_{l}})}^{k-1}a_{s}\mathbbm{1}_{N}^T\epsilon_{s+1}^{(1)}. \label{SumRIkSigNk}
\end{align}
The second term at the right-hand side of \eqref{SumRIkSigNk} tends to zero as $k\to\infty$ since $\|\sum_{k=1}^{\infty}a_{k}\epsilon_{k}^{(1)}\|<\infty$ by A3. We now show that the first term at the right-hand side of \eqref{SumRIkSigNk} converges:
\begin{align}\label{UBarRIKToUstar}
\mathbbm{1}_{N}^{T}\overline{\mathbbm{u}}_{r(i_{l},\overline{\sigma}_{n_{l}})}\xrightarrow[k\to\infty]{} \mathbbm{1}_{N}^{T}\mathbbm{u}^{*}.
\end{align}
Consider the components of $\overline{\mathbbm{u}}_{r(i_{l},\overline{\sigma}_{n_{l}})}$. First, by \eqref{AuxiSigBar} it follows that
\begin{align}\label{UBarIKSigBar}
\overline{u}_{i_{l},r(i_{l},\overline{\sigma}_{n_{l}})} =u_{i_{l},r(i_{l},\overline{\sigma}_{n_{l}})}=u_{i_{l}}^{*},
\end{align}
while for $j\ne i_{l}$ we have $r(\overline{\sigma}_{n_{l}})\le r(j,\overline{\sigma}_{n_{l}})\le r(i_{l},\overline{\sigma}_{n_{l}})< r(\overline{\sigma}_{n_{l}}+1)$, and hence
\begin{align}
\overline{u}_{j,r(i_{l},\overline{\sigma}_{n_{l}})}&=u_{j,r(i_{l},\overline{\sigma}_{n_{l}})} \notag\\ &=u_{j,r(j,\overline{\sigma}_{n_{l}})}+\sum_{s=r(j,\overline{\sigma}_{n_{l}})}^ {r(i_{l},\overline{\sigma}_{n_{l}})-1}a_{s}O_{j,s+1}\notag\\ &=u_{j}^{*}+\sum_{s=r(j,\overline{\sigma}_{n_{l}})}^ {r(i_{l},\overline{\sigma}_{n_{l}})-1}a_{s}O_{j,s+1}. \label{UBarJRIKSigBar}
\end{align}
Noticing $r(i_{l},\overline{\sigma}_{n_{l}})-r(j,\overline{\sigma}_{n_{l}})\le r(i_{l},\overline{\sigma}_{n_{l}})-r(\overline{\sigma}_{n_{l}})\le d(\mathcal{G})$ and \eqref{AkOkToLn}, we conclude that the second term at the right-hand side of \eqref{UBarJRIKSigBar} tends to zero as $k\to\infty$. This combined with \eqref{UBarIKSigBar} implies \eqref{UBarRIKToUstar}. Therefore, by \eqref{SumRIkSigNk} we have $\mathbbm{1}_{N}^{T}\overline{\mathbbm{u}}_{k}\xrightarrow[k\to\infty]{} \mathbbm{1}_{N}^{T}\mathbbm{u}^{*}$ for $k\in(r(i_{l},\overline{\sigma}_{n_{l}}), n_{l+1})$.

Combining 1) and 2) we derive \eqref{DistanceForUkBarSstar0}.

\textbf{Step 2:} Consider $\omega \in \Omega$ where \eqref{SumAkEpijk} holds. Let the sequence $\{\overline{\mathbbm{u}}_{k}\}$ be generated by \eqref{UkBarHUkBar} and \eqref{SigmaBarHUkBar} and let $S$ be a closed subset of $\mathbb{R}^N$ such that $J\cap S\ne\emptyset$. Assume $\overline{\mathbbm{u}}_{k} \in S$ for $k \ge k_{0}$ ʱ. If $[\delta_{1},\delta_{2}]$ is an interval such that $d([\delta_{1},\delta_{2}],v(J \cap S))>0$, then for any bounded subsequence $\{\overline{\mathbbm{u}}_{n_{k}}\}_{n_{k} \ge k_{0}}\subset\{\overline{\mathbbm{u}}_{k}\}_{k\ge k_{0}}$, $\{v(\overline{\mathbbm{u}}_{n})\}_{n \ge k_{0}}$ cannot cross $[\delta_{1},\delta_{2}]$ infinitely many times with starting points $\overline{\mathbbm{u}}_{n_{k}}$ where "$\{v(\overline{\mathbbm{u}}_{n})\}_{n \ge k_{0}}$ crosses $[\delta_{1},\delta_{2}]$ with starting points $\overline{\mathbbm{u}}_{n_{k}}$" means that $v(\overline{\mathbbm{u}}_{n_{k}})\le\delta_{1}$ and there exists $l_{k}>n_{k}$ such that $v(\overline{\mathbbm{u}}_{l_{k}})\ge\delta_{2}$ and $\delta_{1}<v(\overline{\mathbbm{u}}_{n})<\delta_{2}$ for $n: n_{k}<n<l_{k}$.

In what follows $n$ and $n_{l}$ are always assumed to be equal to or greater than $n_{0}$. Assume the converse that $v(\overline{\mathbbm{u}}_{n})$ crosses $[\delta_{1},\delta_{2}]$ infinitely many times with starting points $\overline{\mathbbm{u}}_{n_{k}}$. Without loss of generality, we may assume $\{\overline{\mathbbm{u}}_{n_{k}}\}$ is convergent: $\lim_{k \to \infty}\overline{\mathbbm{u}}_{n_{k}} \triangleq \overline{\mathbbm{u}} \in S$ where $\overline{u}\in S$ because $S$ is closed. From
\eqref{AkOkToLn} and \eqref{OBar} we have
\begin{align}\label{Unk1MunisUnk}
\|\overline{\mathbbm{u}}_{n_{k}+1}-\overline{\mathbbm{u}}_{n_{k}}\| =\|a_{n_{k}}\overline{\mathbbm{O}}_{n_{k}+1}\|
\xrightarrow[k\to\infty]{}0.
\end{align}
By the definition of crossing, $v(\overline{\mathbbm{u}}_{n_{k}})\le\delta_{1}<v(\overline{\mathbbm{u}}_{n_{k}+1})$ which associated with \eqref{Unk1MunisUnk} implies $\lim_{k\to\infty}v(\overline{\mathbbm{u}}_{n_{k}})=\delta_{1}=v(\overline{\mathbbm{u}})$.
In view of $d([\delta_{1},\delta_{2}],v(J \cap S))>0$ and $\overline{\mathbbm{u}}\in S$, we conclude that
\begin{align}\label{DistanceUBarJ}
d(\overline{\mathbbm{u}},J)>0.
\end{align}
This is because if $d(\overline{\mathbbm{u}},J)=0$, then $d(\overline{\mathbbm{u}}_{n_{k}},J)\xrightarrow[k\to\infty]{}0$. However, $d(\overline{\mathbbm{u}}_{n_{k}},S)\xrightarrow[k\to\infty]{}0$, so $d(\overline{\mathbbm{u}}_{n_{k}},J\cap S)\xrightarrow[k\to\infty]{}0$, which implies that $d(v(\overline{\mathbbm{u}}_{n_{k}}),v(J\cap S))\xrightarrow[k\to\infty]{}0$. This contradicts to  $\lim_{k\to\infty}v(\overline{\mathbbm{u}}_{n_{k}})=\delta_{1}$ and
$d([\delta_{1},\delta_{2}],v(J \cap S))>0$.

When $T>0$ is sufficiently small, by Lemma~\ref{PropertyConSubSeq} there is no truncation for \eqref{UkBarHUkBar} and \eqref{SigmaBarHUkBar} with time $k$ running from $n_{k}$ to $m(n_{k},T)$ for sufficiently large $n_{k}$, i.e., $m(n_{k},T)< r(\overline{\sigma}_{n_{k}}+1)$. Assume $k$ is large enough so that
\eqref{UsPlusoneMinusUnk} also holds. Then, we have
\begin{align}\label{VuBarMnkT}
&v(\overline{\mathbbm{u}}_{m(n_{k},T)+1})-v(\overline{\mathbbm{u}}_{n_{k}}) =\sum_{l=n_{k}}^{m(n_{k},T)}a_{l} \overline{\mathbbm{O}}_{l+1}^T \nabla_{\mathbbm{u}}v(\overline{\mathbbm{u}})+o(T) \notag \\
&=\sum_{l=n_{k}}^{m(n_{k},T)}a_{l}\mathbbm{h}^T(\overline{\mathbbm{u}})(P-D) \nabla_{\mathbbm{u}}v(\overline{\mathbbm{u}})
+\sum_{l=n_{k}}^{m(n_{k},T)}a_{l}\overline{\epsilon}_{l+1} \notag \\
&+\sum_{l=n_{k}}^{m(n_{k},T)}a_{l}(\mathbbm{h}(\overline{\mathbbm{u}}_{l}) -\mathbbm{h}(\overline{\mathbbm{u}})) (P-D)\nabla_{\mathbbm{u}}v(\overline{\mathbbm{u}})
+o(T).
\end{align}
By \eqref{EstForUBarsPlusOne} and Lemma~\ref{PropertyEpsiBar} (where $C>M+1$), it follows that
\begin{align}\label{AlEpslnktoMnkT}
\sum_{l=n_{k}}^{m(n_{k},T)}a_{l}\overline{\epsilon}_{l+1}=o(T)
\end{align}
for sufficiently small $T>0$ and large enough $k$. By Lemma~\ref{PropertyConSubSeq}, the continuity of $h(\cdot)$, and \eqref{EstForUBarsPlusOne} we have
\begin{align}\label{AlNablaVUBar}
\sum_{l=n_{k}}^{m(n_{k},T)}a_{l}(\mathbbm{h}(\overline{\mathbbm{u}}_{l}) -\mathbbm{h}(\overline{\mathbbm{u}}))(P-D) \nabla_{\mathbbm{u}}v(\overline{\mathbbm{u}})
=o(T).
\end{align}
This incorporating with \eqref{VuBarMnkT}, \eqref{AlEpslnktoMnkT} and \eqref{AlNablaVUBar} leads to
\begin{align}\label{VumkMinusVunk}
&v(\overline{\mathbbm{u}}_{m(n_{k},T)+1})-v(\overline{\mathbbm{u}}_{n_{k}})\notag\\
&=\sum_{l=n_{k}}^{m(n_{k},T)}a_{l}\mathbbm{h}^T(\overline{\mathbbm{u}})(P-D) \nabla_{\mathbbm{u}}v(\overline{\mathbbm{u}})+o(T).
\end{align}
Thus , by \eqref{NablaVPIHUZeta}, \eqref{DistanceUBarJ}, and \eqref{VumkMinusVunk} there exists a $\zeta>0$ such that
\begin{align}\label{VUBarMnkTMinusVUBarnk}
v(\overline{\mathbbm{u}}_{m(n_{k},T)+1})-v(\overline{\mathbbm{u}}_{n_{k}})<-\frac{\zeta}{2} T
\end{align}
for sufficiently large $k$. On the other hand, by Lemma ~\ref{PropertyConSubSeq} we have
\[
\max_{n_{k} \le s \le m(n_{k},T)}|v(\overline{\mathbbm{u}}_{s+1})-v(\overline{\mathbbm{u}}_{n_{k}})|\xrightarrow[T \to 0]{} 0
\]
which means that $v(\overline{\mathbbm{u}}_{m(n_{k},T)+1})\in [\delta_{1},\delta_{2}]$ for sufficiently small $T>0$. This is a contradiction to \eqref{VUBarMnkTMinusVUBarnk}.

\textbf{Step 3:}
By \eqref{ChoiseForC0} and the definition of $v(\cdot)$ we have
\begin{align*}
&\inf_{\|\mathbbm{u}\|_{\infty}=c_{0}}v(\mathbbm{u})=\min_{i \in \mathcal{N}}\Big\{\min \{ \int_{u_{i}^0}^{c_{0}}h_{i}(t)\text{d}t ,\notag \\
&\int_{u_{i}^0}^{-c_{0}}h_{i}(t) \text{d}t\}\Big\} > \sum_{j \in \mathcal{N}}\int_{u_{j}^0}^{u_{j}^*}h_{j}(t) \text{d}t \notag \\
&=v(\mathbbm{u}^*).
\end{align*}
By Lemma~\ref{SingleRoot}, $J \cap S^*(0)$ is a singleton. Thus, there exists a nonempty interval $[\delta_{1}, \delta_{2}]$ such that $d([\delta_{1},\delta_{2}],v(J \cap S^*(0)))>0$ and $[\delta_{1},\delta_{2}] \subset (v(\mathbbm{u}^*),\inf_{\|\mathbbm{u}\|_{\infty}=c_{0}}v(\mathbbm{u}))$. So, there exists a sufficiently small $\epsilon >0$ such that $d([\delta_{1},\delta_{2}],v(J \cap S^*(\epsilon)))>0$. By \eqref{DistanceForUkBarSstar0}, there exists $k_{0}$ such that $\overline{\mathbbm{u}}_{k} \in S^*(\epsilon)$ for $k \ge k_{0}$. Since $M_{k} \xrightarrow[k \to \infty]{} \infty$, there exists $k_{1}$ such that $M_{\overline{\sigma}_{n_{k}}}>c_{0}$ for $k > k_{1}$. By the definition of $\{n_{k}\}$ we have
\[
\|\overline{\mathbbm{u}}_{n_{k+1}-1}+a_{n_{k+1}-1} \overline{\mathbbm{O}}_{n_{k+1}}\|_{\infty} >M_{\overline{\sigma}_{n_{k}}}
\]
for $k > k_{0} \vee k_{1}$. This means that $\{v(\overline{\mathbbm{u}}_{n}), n > k_{0} \vee k_{1}\}$ crosses $[\delta_{1},\delta_{2}]$ infinitely many times with starting points $\overline{\mathbbm{u}}_{n_{k}}=\mathbbm{u}^{*}$, but it is impossible by assertion in \textbf{Step 2}. Thus, the lemma is proved. \qed
\end{pf}

Set $\overline{\sigma}=\overline{\sigma}_{K}$ where $K$ is the smallest time when the algorithm defined by \eqref{UkBarHUkBar} and \eqref{SigmaBarHUkBar} has no more truncations. For the samples $\omega \in \Omega$ where \eqref{SumAkEpijk} holds, Lemma~\ref{FiniteTruncation} means that there are only a finite number of truncations for \eqref{UkBarHUkBar} and \eqref{SigmaBarHUkBar}. At the same time, noticing \eqref{AuxiUBar} and \eqref{AuxiSigBar}, we conclude that there are only finitely number of times for which $\{\overline{\mathbbm{u}}_{k}\}$ may differ from $\{\mathbbm{u}_{k}\}$, and $\overline{\mathbbm{u}}_{k}=\mathbbm{u}_{k}$ for $k\ge K$. For $\epsilon \ge 0$ set $S(\epsilon)\triangleq\{\mathbbm{u}:\mathbbm{1}_{N}^T \mathbbm{u} \in [\mathbbm{1}_{N}^T \mathbbm{u}^*+\sum_{s=K}^{\infty}a_{s}\mathbbm{1}_{N}^T\overline{\epsilon}_{s+1}-\epsilon, \mathbbm{1}_{N}^T\mathbbm{u}^*+\sum_{s=K}^{\infty}a_{s}\mathbbm{1}_{N}^T\overline{\epsilon}_{s+1} +\epsilon]\}$. After the truncation having ceased the algorithm becomes
\begin{align*}
\overline{\mathbbm{u}}_{k+1}=\overline{\mathbbm{u}}_{k}+a_{k}(-L\mathbbm{y}_{k+1}+\epsilon_{k+1}^{(1)}), \text{ }k\ge K,
\end{align*}
where $\epsilon_{k+1}^{(1)}=[\epsilon_{1,k+1}^{(1)},\cdots,\epsilon_{N,k+1}^{(1)}]$. Therefore
\[
\sum_{s=K}^{\infty}a_{s}\mathbbm{1}_{N}^T\overline{\epsilon}_{s+1}=\sum_{s=K}^{\infty}a_{s} \mathbbm{1}_{N}^T\epsilon_{s+1}^{(1)}< \infty,
\]
and hence $S(\epsilon)$ is well defined. Similar to \textbf{Step 1} in the proof of Lemma~\ref{FiniteTruncation}, we have $d(\overline{\mathbbm{u}}_{k},S(0))\xrightarrow[k \to \infty]{} 0$.

\begin{thm}\label{MainResult}
Assume A1-A4 hold. Then applying the algorithm \eqref{SigmakHat}-\eqref{AlgorithmSigmak} to the multi-agent systems composed of the Hammerstein systems \eqref{HSystem} and the Wiener systems \eqref{WSystem} leads to consensus:
\begin{align}\label{ConsensusForOutput}
y_{i,k} \xrightarrow[k \to \infty]{} y^0 \text{ a.s.}, \forall i \in \mathcal{N},
\end{align}
where $y^0=y^0(\omega)$ may depend on samples $\omega$ and is such that $|y^0|<\infty \text{ a.s.}$, and $ y^0\mathbbm{1}_{N}=\mathbbm{h}(\mathbbm{u}^{0})$, where
\[
h_{i}(u)=\begin{cases}\frac{d_{i}}{c_{i}}f_{i}(u), & \text{if }i\text{ is the Hammerstein system},\\
f_{i}(\frac{d_{i}}{c_{i}}u), & \text{if }i\text{ is the Wiener system}.\end{cases}
\]
\end{thm}
\begin{pf}
We first note that $J\cap S(0)\triangleq\{\mathbbm{u}^{0}\}$ is a singleton by Lemma~\ref{SingleRoot}. By Lemma~1 in \cite{Chen2007}, if $\lim_{k\to\infty}u_{i,k}=u_{i}$, then $\lim_{k\to\infty}y_{i,k}=h_{i}(u_{i})$. Therefore, to prove \eqref{ConsensusForOutput} it suffices to show that $d(\mathbbm{u}_{k},J \cap S(0)) \xrightarrow[k \to \infty]{} 0$ or $\mathbbm{u}_{k}\xrightarrow[k \to \infty]{} \mathbbm{u}^{0}$. From Lemma~\ref{FiniteTruncation} we know that $\{\overline{\mathbbm{u}}_{k}\}$ is bounded. Set
\begin{align*}
v_{1}\triangleq \liminf_{k \to \infty}v(\overline{\mathbbm{u}}_{k}) \le \limsup_{k \to \infty}v(\overline{\mathbbm{u}}_{k}) \triangleq v_{2}.
\end{align*}

If $v_{1}<v_{2}$, then there exist $\delta_{1}$ and $\delta_{2}$ such that $\delta_{1}<\delta_{2}$, $d([\delta_{1},\delta_{2}],v(J\cap S(0))>0$, and $[\delta_{1},\delta_{2}] \subset (v_{1},v_{2})$. Thus, there exist $\epsilon>0$ and $k_{0}$ such that $d([\delta_{1},\delta_{2}],v(J \cap S(\epsilon)))>0$ and $\overline{\mathbbm{u}}_{k} \in S(\epsilon),k>k_{0}$. This implies that $\{v(\overline{\mathbbm{u}}_{k}),k > k_{0}\}$ crosses $[\delta_{1},\delta_{2}]$ infinitely many times and contradicts to what proved in \textbf{Step 2} in the proof of Lemma~\ref{FiniteTruncation}. So, $v_{1}=v_{2}$, i.e., $\{v(\overline{\mathbbm{u}}_{k})\}$ is convergent.

Assume the converse: there exists $\{\overline{\mathbbm{u}}_{n_{k}}\} \subset \{\overline{\mathbbm{u}}_{k}\}$ such that $\lim_{k \to \infty}\overline{\mathbbm{u}}_{n_{k}}=\overline{\mathbbm{u}} \ne \mathbbm{u}^{0}$. Noticing $d(\overline{\mathbbm{u}}_{k},S(0))\xrightarrow[k \to \infty]{} 0$, similar to
\eqref{VUBarMnkTMinusVUBarnk} we have $v(\overline{\mathbbm{u}}_{m(n_{k},T)+1})-v(\overline{\mathbbm{u}}_{n_{k}}) \le -\frac{\zeta}{2}T$ for sufficiently small $T>0$ and large enough $k$. This contradicts to the convergence of $v(\overline{\mathbbm{u}}_{k})$. Therefore, $d(\overline{\mathbbm{u}}_{k},J \cap S(0)) \xrightarrow[k \to \infty]{} 0$. By Lemma~\ref{FiniteTruncation}, $\{\mathbbm{u}_{k}\}$ may differ from $\{\overline{\mathbbm{u}}_{k}\}$ only by a finite number of terms, so $d(\mathbbm{u}_{k},J \cap S(0)) \xrightarrow[k \to \infty]{} 0$. \qed
\end{pf}

\section{Numerical simulation}
Consider the undirected communication graph with four agents presented in Fig.~1. The corresponding Laplace matrix is as follows:
\[
L=\begin{bmatrix} 2 & -1 & 0 & -1 \\
-1 & 3 & -1 & -1 \\
0 & -1 & 1 & 0 \\
-1 & -1 & 0 & 2
\end{bmatrix}.
\]
\begin{figure}
\begin{center}
\includegraphics[height=4cm]{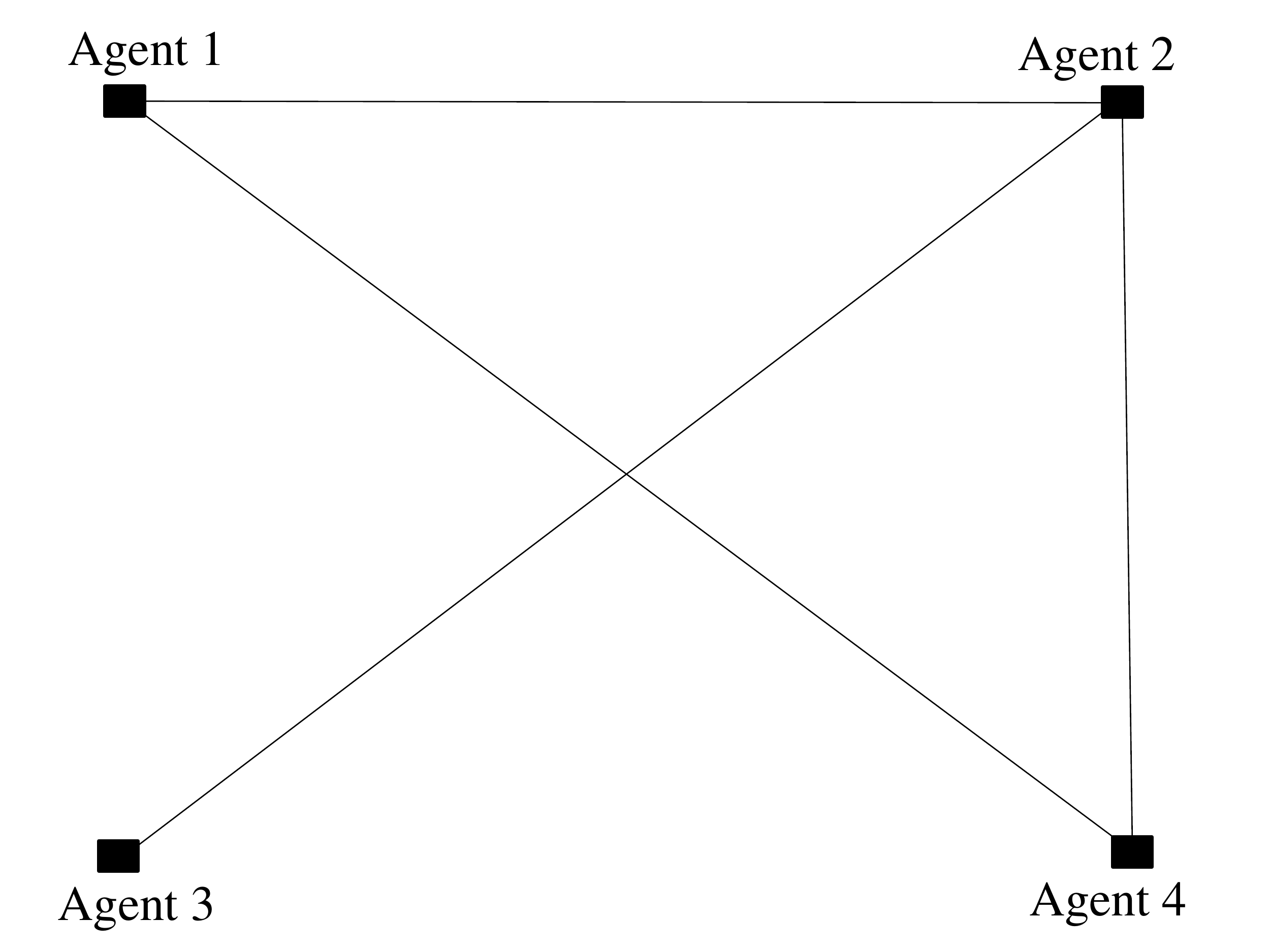}    
\caption{Simulation graph}  
\label{fig1}                                 
\end{center}                                 
\end{figure}
Using the algorithm given by \eqref{SigmakHat}-\eqref{AlgorithmSigmak}, we compute the output of the systems for three cases: 1. All four agents are Hammerstein systems (H) but with different parameters and static functions. 2. All four agents are Wiener systems (W) but with different parameters and static functions. 3. Two are Hammerstein systems and the other two are Wiener systems. To be precise, they are as follows.

Case 1.
\begin{align*}
H:&\begin{cases}v_{1,k}=f_{1}(u_{1,k})=-u_{1,k}^3-u_{1,k}, \\ y_{1,k+1}+0.2y_{1,k}+0y_{1,k-1}+0.6y_{1,k-2}=v_{1,k}\\-0.3v_{1,k-1}-1.2v_{1,k-2};\end{cases}\\
H:&\begin{cases}v_{2,k}=f_{2}(u_{2,k})=-2u_{2,k}+1, \\
y_{2,k+1}+0.6y_{2,k}+0.5y_{2,k-1}+0.4y_{2,k-2}=v_{2,k}\\-v_{2,k-1}-2v_{2,k-2};\end{cases}\\
H:&\begin{cases}v_{3,k}=f_{3}(u_{3,k})=(u_{3,k}-1)^3, \\ y_{3,k+1}-0.15y_{3,k}+0y_{3,k-1}+0.5y_{3,k-2}=v_{3,k}\\+0.2v_{3,k-1}-0.4v_{3,k-2};\end{cases}\\
H:&\begin{cases}v_{4,k}=f_{4}(u_{4,k})=u_{4,k}^3+1, \\
y_{4,k+1}+0.76y_{4,k}+0.5y_{4,k-1}+0.6y_{4,k-2}=v_{4,k}\\+0.5v_{4,k-1}.\end{cases}
\end{align*}

Case 2.
\begin{align*}
W:&\begin{cases}v_{1,k+1}+0.2v_{1,k}+0v_{1,k-1}+0.6v_{1,k-2}=u_{1,k}\\-0.3u_{1,k-1}-1.2u_{1,k-2}, \\ y_{1,k+1}=f_{1}(v_{1,k+1})=-v_{1,k+1}^3-v_{1,k+1};\end{cases}\\
W:&\begin{cases}v_{2,k+1}+0.6v_{2,k}+0.5v_{2,k-1}+0.4v_{2,k-2}=u_{2,k}\\-u_{2,k-1}-2u_{2,k-2}, \\
y_{2,k+1}=f_{2}(v_{2,k+1})=-2v_{2,k+1}+1;\end{cases}\\
W:&\begin{cases}v_{3,k+1}-0.15v_{3,k}+0v_{3,k-1}+0.5v_{3,k-2}=u_{3,k}\\+0.2u_{3,k-1}-0.4u_{3,k-2}, \\ y_{3,k+1}=f_{3}(v_{3,k+1})=(v_{3,k+1}-1)^3;\end{cases}\\
W:&\begin{cases}v_{4,k+1}+0.76v_{4,k}+0.5v_{4,k-1}+0.6v_{4,k-2}=u_{4,k}\\+0.5u_{4,k-1}, \\
y_{4,k+1}=f_{4}(v_{4,k+1})=v_{4,k+1}^3+1.\end{cases}
\end{align*}

Case 3.
\begin{align*}
W:&\begin{cases}v_{1,k+1}+0.2v_{1,k}+0v_{1,k-1}+0.6v_{1,k-2}=u_{1,k}\\-0.3u_{1,k-1}-1.2u_{1,k-2}, \\ y_{1,k+1}=f_{1}(v_{1,k+1})=-v_{1,k+1}^3-v_{1,k+1};\end{cases}\\
W:&\begin{cases}v_{2,k+1}+0.6v_{2,k}+0.5v_{2,k-1}+0.4v_{2,k-2}=u_{2,k}\\-u_{2,k-1}-2u_{2,k-2}, \\
y_{2,k+1}=f_{2}(v_{2,k+1})=-2v_{2,k+1}+1;\end{cases}\\
H:&\begin{cases}v_{3,k}=f_{3}(u_{3,k})=(u_{3,k}-1)^3, \\ y_{3,k+1}-0.15y_{3,k}+0y_{3,k-1}+0.5y_{3,k-2}=v_{3,k}\\+0.2v_{3,k-1}-0.4v_{3,k-2};\end{cases}\\
H:&\begin{cases}v_{4,k}=f_{4}(u_{4,k})=u_{4,k}^3+1, \\
y_{4,k+1}+0.76y_{4,k}+0.5y_{4,k-1}+0.6y_{4,k-2}=v_{4,k}\\+0.5v_{4,k-1}.\end{cases}
\end{align*}

All observation noises $\{\epsilon_{12,k}\}, \{\epsilon_{14,k}\}, \{\epsilon_{21,k}\}, \{\epsilon_{23,k}\}$, $ \{\epsilon_{24,k}\}, \{\epsilon_{32,k}\}, \{\epsilon_{41,k}\}$, and $\{\epsilon_{42,k}\}$ are mutually independent and normally distributed with with zero mean and variance $1$. It can straightforwardly be verified that A1-A4 hold. We take all initial values to be $0$ and $u_{1}^{*}=1, u_{2}^{*}=2, u_{3}^{*}=3, u_{4}^{*}=4$, $c_{M}=55$. The algorithm \eqref{SigmakHat} - \eqref{AlgorithmSigmak} is applied. Simulation results are presented in Fig.~2 - Fig.~7 for cases~1-3, respectively. From the figures it is seen that the output consensus is achieved and the input at all agents are convergent for all these cases.

\begin{figure}
\begin{center}
\includegraphics[height=4cm]{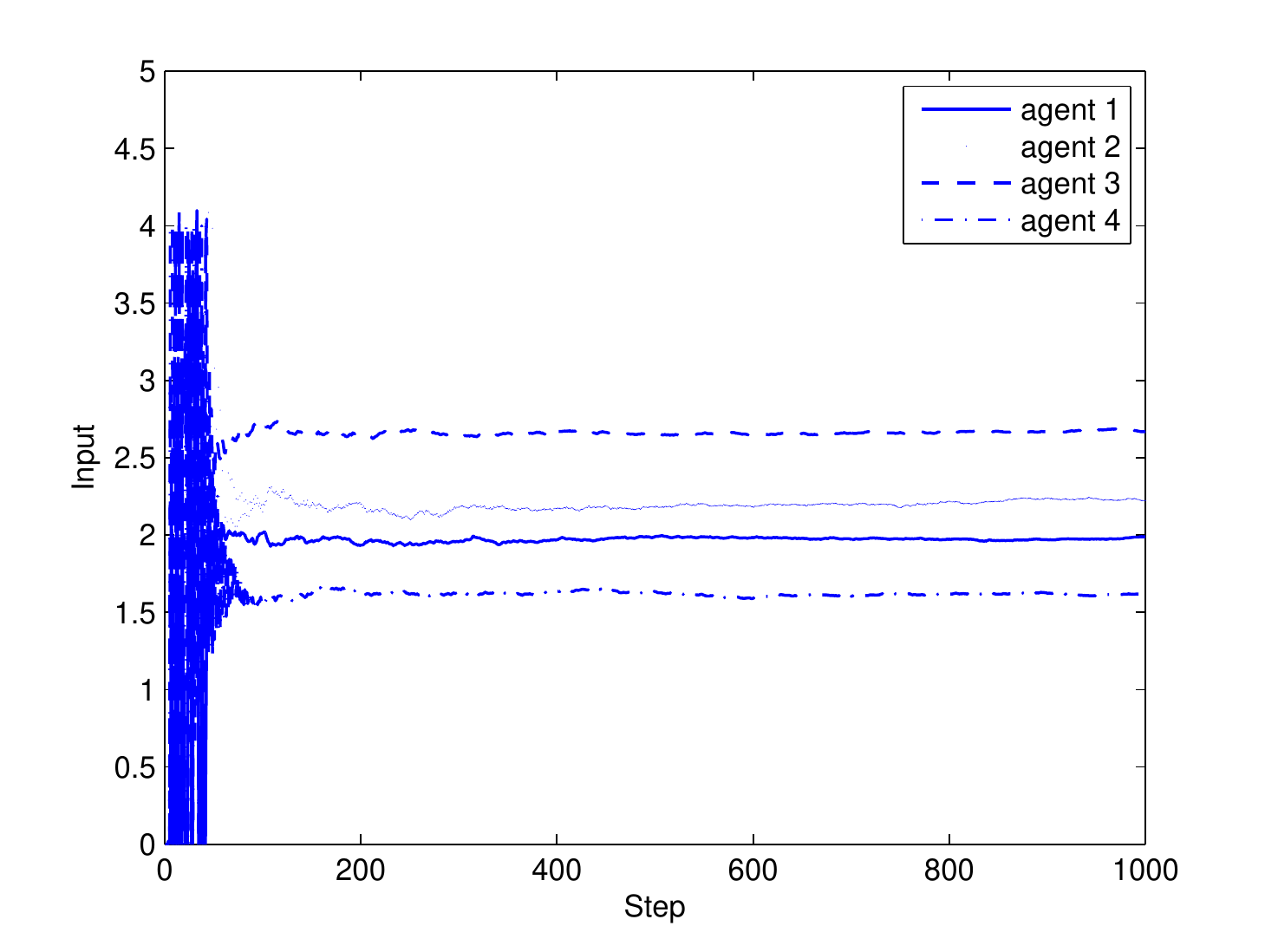}    
\caption{Inputs for Case 1}  
\label{fig1}                                 
\end{center}                                 
\end{figure}

\begin{figure}
\begin{center}
\includegraphics[height=4cm]{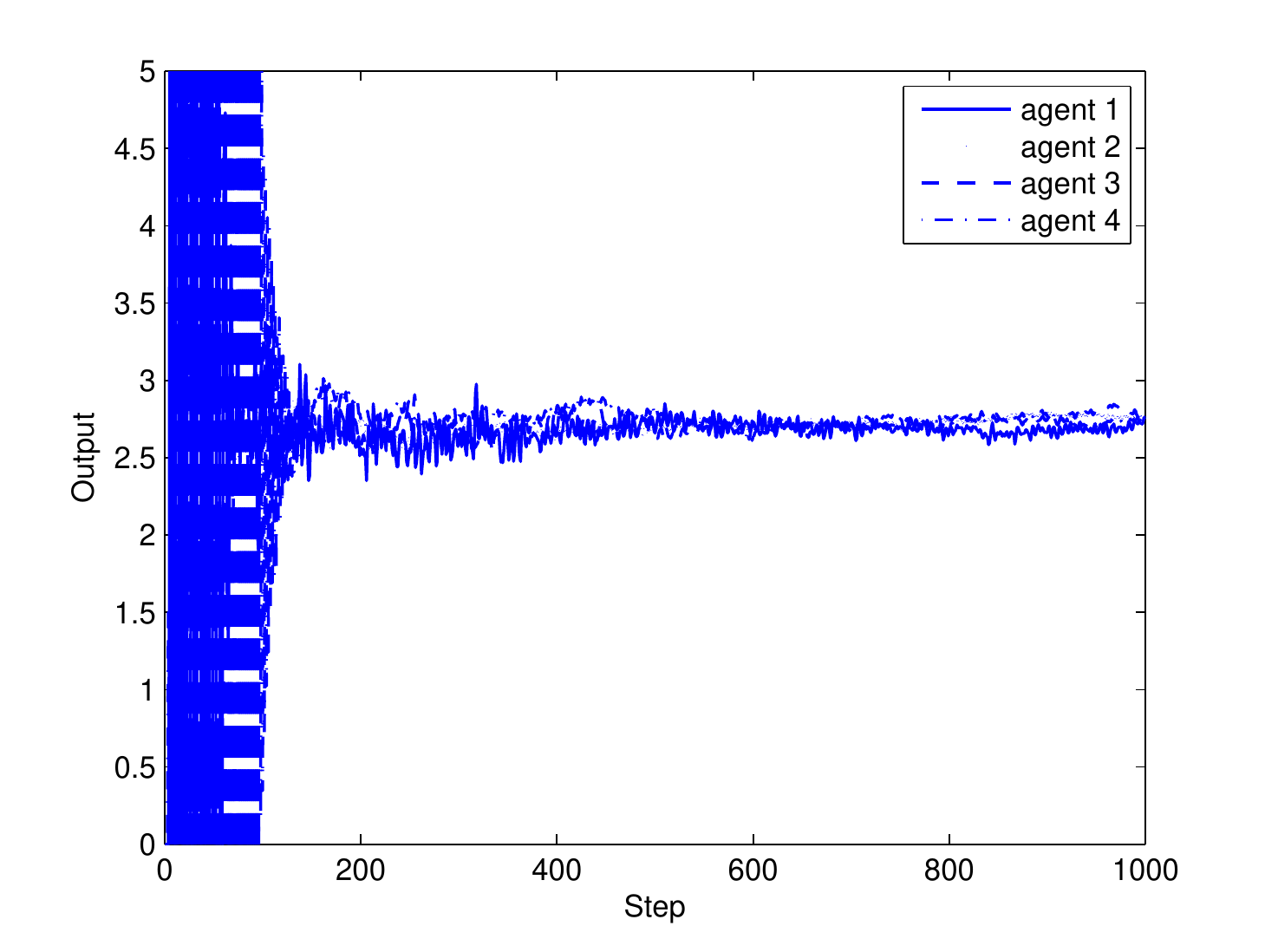}    
\caption{Outputs for Case 1}  
\label{fig1}                                 
\end{center}                                 
\end{figure}

\begin{figure}
\begin{center}
\includegraphics[height=4cm]{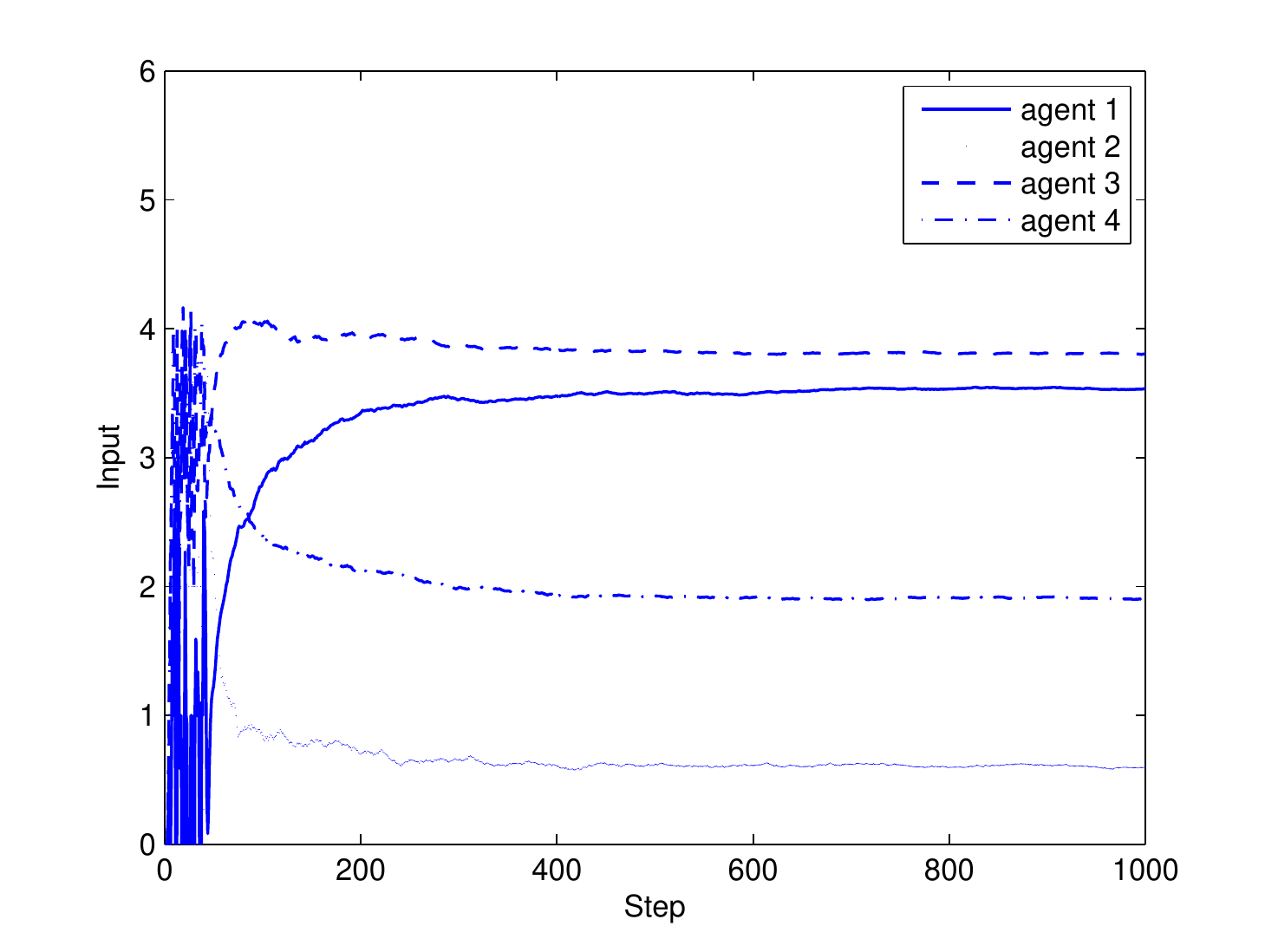}    
\caption{Inputs for Case 2}  
\label{fig1}                                 
\end{center}                                 
\end{figure}

\begin{figure}
\begin{center}
\includegraphics[height=4cm]{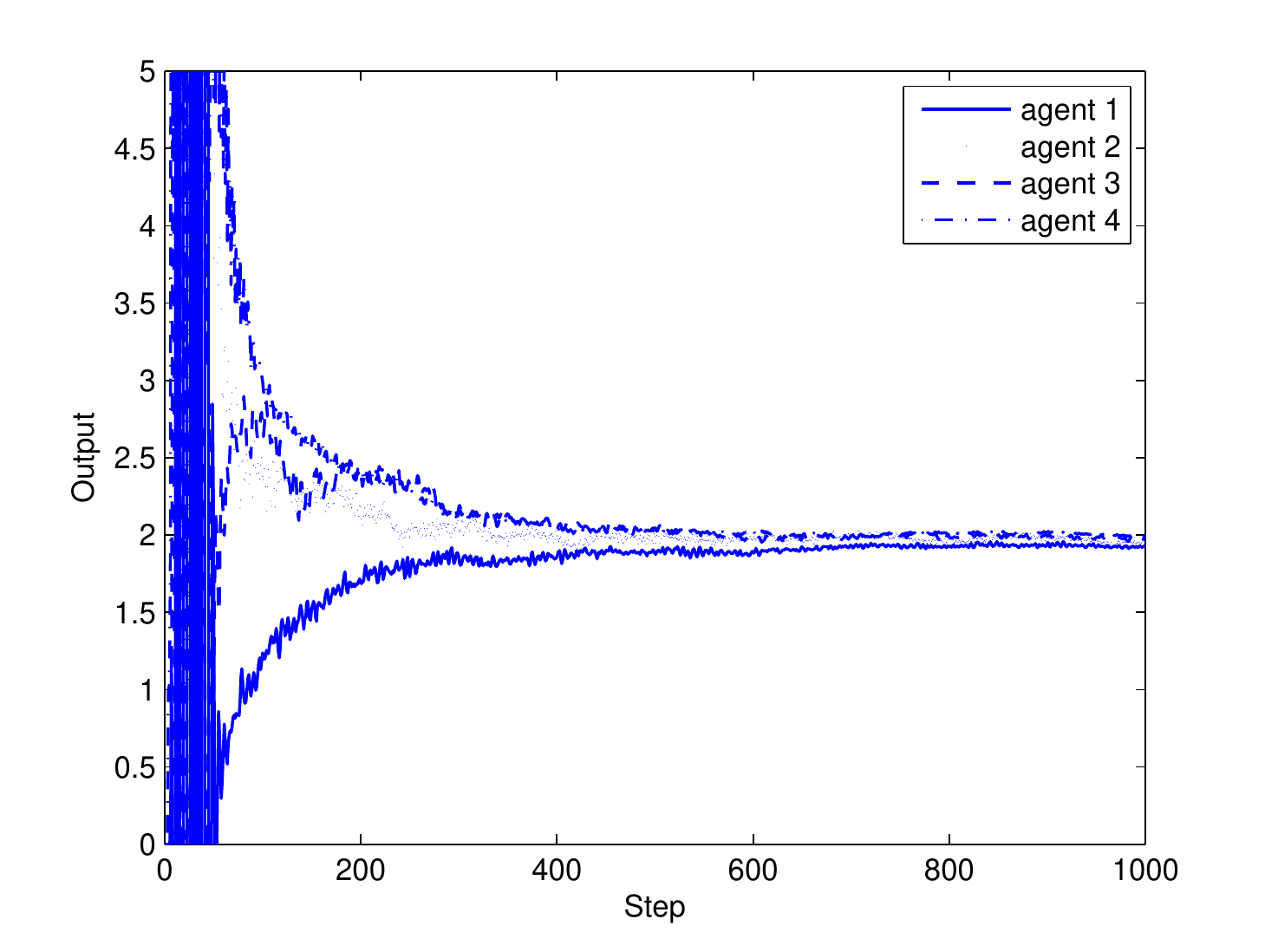}    
\caption{Outputs for Case 2}  
\label{fig1}                                 
\end{center}                                 
\end{figure}

\begin{figure}
\begin{center}
\includegraphics[height=4cm]{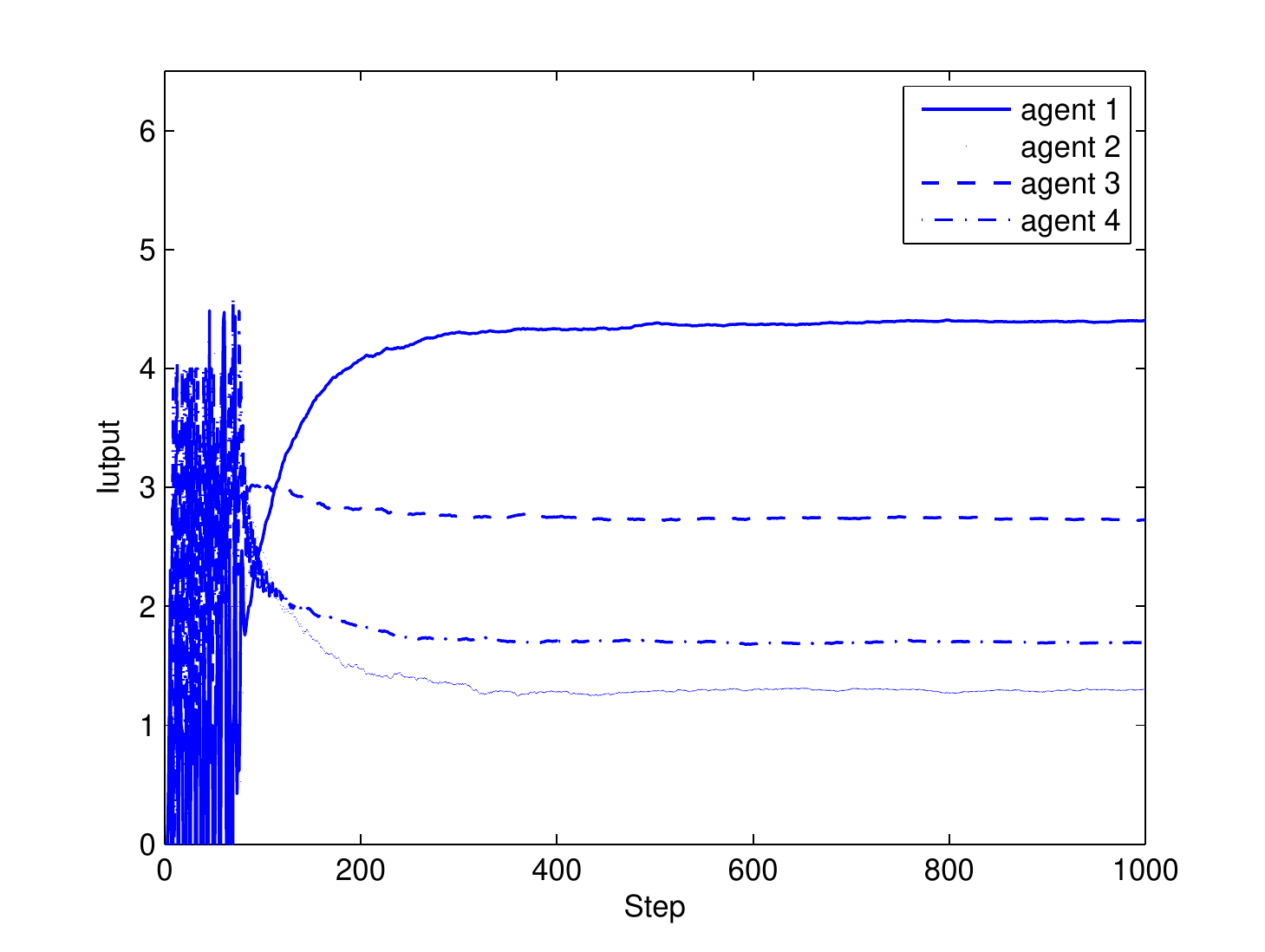}    
\caption{Inputs for Case 3}  
\label{fig1}                                 
\end{center}                                 
\end{figure}

\begin{figure}
\begin{center}
\includegraphics[height=4cm]{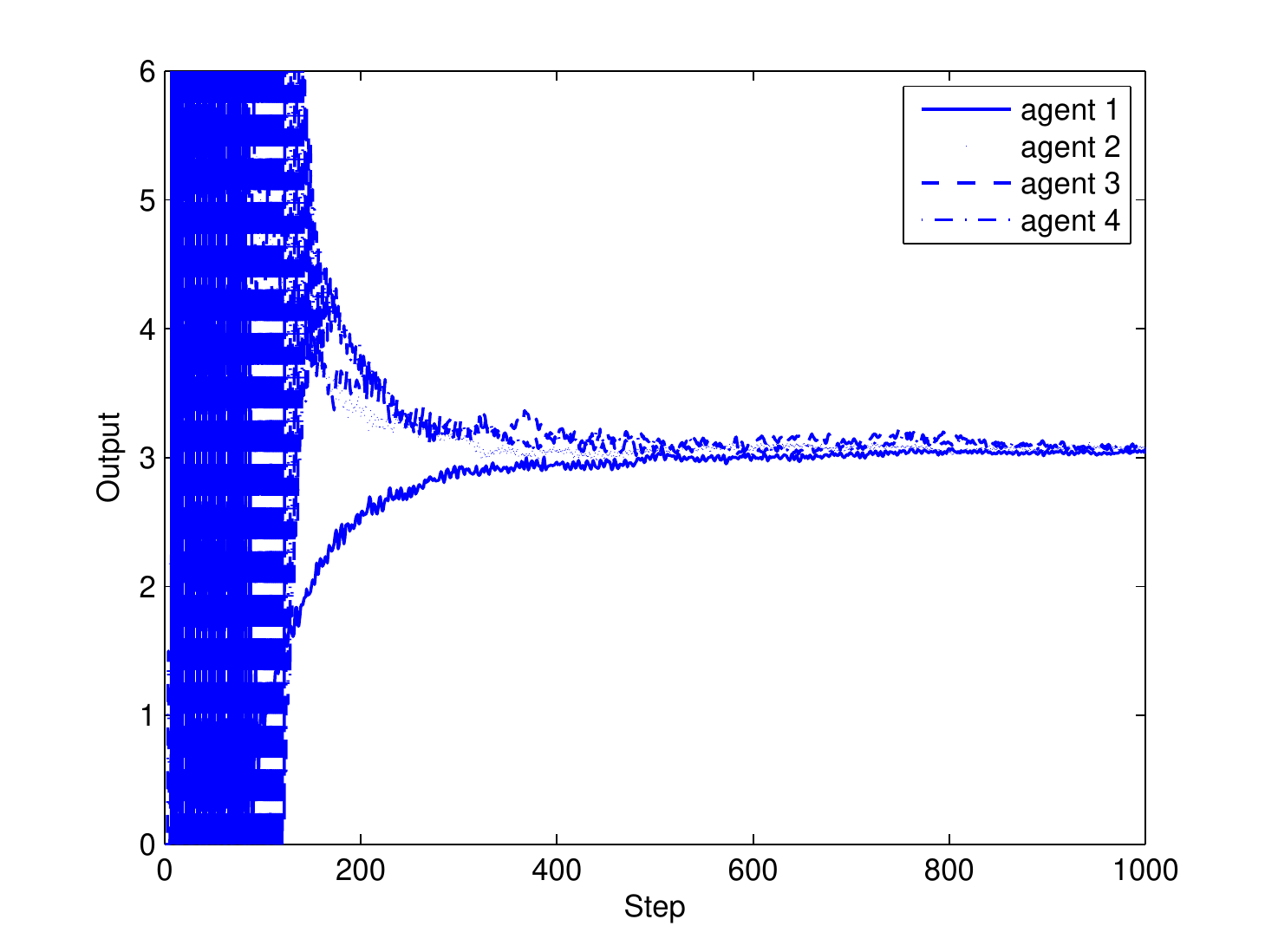}    
\caption{Outputs for Case 3}  
\label{fig1}                                 
\end{center}                                 
\end{figure}


\section{Conclusions}
In this paper the output consensus problem of networked Hammerstein and Wiener systems is studied in a noisy communication environment. Each agent being a Hammerstein or Wiener system is assumed to be open-loop stable, and its static nonlinearity is allowed to grow up but not faster than a polynomial. A control algorithm based on DSAAWET is proposed. The algorithm is transformed to a centralized SAAWET by introducing auxiliary sequences, but its convergence is established by a treatment different from a direct application of the convergence theorem presented in \cite{Chen2002}. It is shown that the proposed algorithm leads to the output consensus with probability one. For further research it is of interest to consider the output consensus problem for other subsystems, for example, the nonlinear ARX systems, the multi-input multi-output systems, etc. It is also of interest to consider other kind of communication graphs.

\begin{ack}                               
The authors would like to thank Dr. W. X. Zhao for his useful suggestions and kind help. 
\end{ack}


\small
\begin{thebibliography}{99}     
\bibitem[Jadbabaie {\em et al.} (2003)]{Jadbabaie2003}
Jadbabaie, A., Lin, J., \& Morse, A. S. (2003). Coordination of groups of mobile autonomous agents using nearest neighbor rules. {\it IEEE Trans. Autom. Control} 48(6), 988-1001.

\bibitem[Vicsek {\em et al.} (1995)]{Vicsek1995}
Vicsek, T., Czirk, Benjacob, A. E., Cohen, I. I. \& Shochet, O. (1995). Novel type of phase transition in a system of self-driven particles. {\it Physical Review Letters,} 75(6), 1226.

\bibitem[Saber {\em et al.} (2004)]{Saber2004}
Olfati-Saber, R., Fax, J. A. \& Murray, R. M. (2004). Consensus problems in networks of agents with switching topology and time-delays. {\it IEEE Trans. Autom. Control,} 49(9), 1520-1533.

\bibitem[Ren and Beard (2005)]{Ren2005}
Ren, W., \& Beard, R. W. (2005). Consensus seeking in multiagent systems under dynamically changing interaction topologies. {\it IEEE Trans. Autom. Control,} 50, 655-661.


\bibitem[Li and Zhang (2009)]{Li2009}
Li, T., \& Zhang, J. F. (2009). Mean square average-consensus under measurement noises and fixed topologies: Necessary and sufficient conditions. {\it Automatica,} 45(8), 1929-1936.

\bibitem[Li and Zhang (2010)]{Li2010}
Li, T., \& Zhang, J. F. (2010). Consensus conditions of multi-agent systems with time-varying topologies and stochastic communication noises. {\it IEEE Trans. Autom. Control,} 55(9), 2043-2057.

\bibitem[Huang and Manton (2009)]{Huang2009}
Huang, M. Y., \& Manton, J. H. (2009). Coordination and consensus of networked agents with noisy measurements: Stochastic algorithms and asymptotic behavior. {\it SIAM Journal on Control and Optimization,} 48(1), 134-161.

\bibitem[Fang {\em et al.} (2012)]{Fang2012}
Fang, H. T., Chen, H.-F., \& Wen, L. (2012). On control of strong consensus for networked agents with noisy observations. {\it Journal of Systems Science Complexity,} 25(1), 1-12.

\bibitem[Ren (2008)]{Ren2008}
Ren, W. (2008). On consensus algorithms for double-integrator dynamics. {\it IEEE Trans. on Autom. Control,} 53(6), 1503-1509.

\bibitem[Chen {\em et al.} (2013)]{ChenY2013}
Chen, Y., Lv, J. H., Yu, X. H., \& Lin, Z. L. (2013). Consensus of discrete-time second-order multiagent systems based on infinite products of general stochastic matrices. {\it SIAM Journal on Control and Optimization,} 51(4), 3274-3301.

\bibitem[Yu {\em et al.} (2017)]{Yu2017}
Yu, W. W., Li, Y., Wen, G. H., Yu, X. H., \& Cao, J. D. (2017). Observer design for tracking consensus in second-order multi-agent systems: Fractional order less than two. {\it IEEE Trans. on Autom. Control,} 62(2), 894-900.

\bibitem[Saber {\em et al.} (2007)]{Saber2007}
Olfati-Saber, R., Fax, J. A., \& Murray, R. M. (2007). Consensus and cooperation in networked multi-agent systems. {\it Proceedings of the IEEE,} 95(1), 215-233.

\bibitem[Hua {\em et al.} (2016)]{Hua2016}
Hua, C.-C., You, X., \& Guan, X.-P. (2016). Leader-following consensus for a class of high-order nonlinear multi-agent systems. {\it Automatica,} 73, 138-144.

\bibitem[Li {\em et al.} (2013)]{Li2013}
Li, Z. K., Ren, W., Liu, X. D., \& Fu, M. Y. (2013). Consensus of multi-agent systems with general linear and lipschitz nonlinear dynamics using distributed adaptive protocols. {\it IEEE Trans. Autom. Control,} 58(7), 1786-1791.

\bibitem[Liu {\em et al.}, 2013]{Liu2013}
Liu, K. E., Xie, G. M., Ren, W., \& Wang, L. (2013). Consensus for multi-agent systems with inherent nonlinear dynamics under directed topologies. {\it Systems Control Letters,} 62(2), 152-162.

\bibitem[Liu and Huang (2017)]{Liu2017}
Liu, W., \& Huang, J. (2017). Adaptive leader-following consensus for a class of higher-order nonlinear multi-agent systems with directed switching networks. {\it Automatica,} 79, 84-92.

\bibitem[Wang {\em et al.} (2017)]{Wang2017}
Wang, W., Wen, C. Y., \& Huang, J. S. (2017). Distributed adaptive asymptotically consensus tracking control of nonlinear multi-agent systems with unknown parameters and uncertain disturbances. {\it Automatica,} 77, 133-142.

\bibitem[Munz {\em et al.} (2011)]{Munz2011}
Munz, U., Papachristodoulou, A., \& Allgower, F. (2011). Robust consensus controller design for nonlinear relative degree two multi-agent systems with communication constraints. {\it IEEE Trans. Autom. Control,} 56(1), 145-151.

\bibitem[Lei (2016)]{Lei2016}
Lei, J. L. (2016). Cooperative Estimation and Optimization over Multi-Agent Networks with Uncertainties. {\it PhD Thesis, The University of Chineses Academy of Sciences.}

\bibitem[Lei and Chen (2015)]{Lei2015}
Lei, J. L., \& Chen, H.-F. (2015). Distributed stochastic approximation algorithm with expanding truncations: Algorithm and applications. {\it arXiv:1410.7180v4.}

\bibitem[Chen (2002)]{Chen2002}
Chen, H.-F. (2002). Stochastic Approximation and Its Applications. Kluwer Academic Publisher, {\it Dordrecht,} The Netherlands.

\bibitem[Chen (2007)]{Chen2007}
Chen, H.-F. (2007). Adaptive regulator for Hammerstein and Wiener systems with noisy observations. {\it IEEE Trans. Autom. Control,} 52, 703-709.

\bibitem[Fang and Chen (2001)]{Fang2001}
Fang, H. T., \& Chen, H.-F. (2001). Asymptotic behavior of asynchronous stochastic approximation. {\it Science in China,} 44(4), 249-258.

\bibitem[Robbins and Monro (1951)]{Robbins1951}
Robbins, H., \& Monro, S. (1951). A stochastic approximation method. {\it Ann. Math. Statist.,} 22, 400-407.
\end{thebibliography}

\appendix
\end{document}